\documentclass[11pt]{article}
\usepackage{graphicx}
\usepackage[a4paper,pdftex=true]{geometry}
\usepackage[latin1]{inputenc}
\usepackage[T1]{fontenc}
\usepackage{aeguill}
\usepackage{amssymb}
\usepackage{amssymb,array}
\usepackage{syntax}
\usepackage{mdwtab}
\usepackage{mathenv}
\usepackage{float}
\usepackage{lscape}
\usepackage{here}
\usepackage[a4paper]{geometry}
\usepackage{amssymb}
\usepackage{bm}
\usepackage{fancyhdr}
\usepackage{graphicx}
\usepackage{natbib}
\usepackage{fancyvrb}

\newtheorem{theorem}{Theorem}
\newtheorem{lemma}[theorem]{Lemma}
\newtheorem{corollary}[theorem]{Corollary}

\newtheorem{remark}{Remark}

\newcommand{\cdf}[2]{\ensuremath{F_{#1}(#2)}}

\newcommand{\cdfEst}[3]{\ensuremath{	\widehat{F}_{#1} \left(#2 ; \Vect{#3} \right)}}
\newcommand{\cdfEstSimp}[1]{\ensuremath{	\widehat{F} \left(#1 \right)}}
\newcommand{\qtEstSimp}[1]{\ensuremath{	\widehat{\xi} \left(#1 \right)}}

\newcommand{\qt}[2]{\ensuremath{\xi_{#1}(#2)}}

\newcommand{\lqt}[3]{\ensuremath{\xi_{#1} \left( \Vect{#2},\Vect{#3}\right)}}

\newcommand{\qtEst}[3]{\ensuremath{	\widehat{\xi}_{#1} \left(#2 ; \Vect{#3} \right)}}

\newcommand{\lqtEst}[3]{\ensuremath{\widehat{\xi} \left( \Vect{#1},\Vect{#2} ; \Vect{#3}\right)}}

\newcommand{\tmEst}[1]{\ensuremath{\overline{\Vect{#1}}^{(\Vect{\beta})}}}
\newcommand{\tm}[1]{\ensuremath{\overline{{#1}}^{(\Vect{\beta})}}}

\newcommand{\pdf}[2]{\ensuremath{f_{#1}(#2)}}

\newcommand{\tauB}{\overline{\tau}}

\newcommand{\oo}[1]{\ensuremath{o\left( #1 \right)}}
\newcommand{\oas}[1]{\ensuremath{o_{a.s.}\left( #1 \right)}}

\newcommand{\OO}[1]{\ensuremath{\mathcal{O}\left( #1 \right)}}
\newcommand{\Oas}[1]{\ensuremath{\mathcal{O}_{a.s.}\left( #1 \right)}}

\renewcommand{\Citet}[1]{\citeauthor{#1} \citeyearpar{#1}}

\newcommand{\Vect}[1]{\bm{#1}}
\newcommand{\RR}{\mathbb{R}}
\newcommand{\ZZ}{\mathbb{Z}}
\newcommand{\Esp}{\mathbf{E}}
\newcommand{\tr}[1]{{#1}^T}
\newcommand{\Mat}[1]{\underline{\Vect{#1}}}
\newcommand{\Prob}[1]{\mathbb{P}\left( #1\right)}
\title{\Large \bf \sc Hurst exponent estimation of locally self-similar Gaussian processes using sample quantiles}

\author{ {\sc By Jean-Fran\c{c}ois Coeurjolly$^1$} \\
{\it University of Grenoble 2, France}}
\date{}
\begin{document}
\maketitle

\begin{center}
{\small \begin{minipage}{12cm}
This paper is devoted to the introduction of a new class of consistent estimators of the fractal dimension of locally self-similar Gaussian processes. These estimators are based on convex combinations of sample quantiles of discrete variations of a sample path over a discrete grid of the interval $[0,1]$. We derive the almost sure convergence and the asymptotic normality for these estimators. The key-ingredient is a Bahadur representation for sample quantiles of non-linear functions of Gaussians sequences with correlation function decreasing as $k^{-\alpha}L(k)$ for some $\alpha>0$ and some slowly varying function $L(\cdot)$. %The effectiveness of our procedure is investigated in a simulation study. We also propose a short simulation to illustrate that such estimators are more robust to additive outliers thant classical ones.
\end{minipage}}
\end{center}

\vspace*{.5cm}

\footnotetext[1]{Supported by a grant from IMAG Project AMOA. \\ {\it AMS 2000 subject classifications:} Primary 60G18; secondary 62G30. \\ {\it Key words and phrases}: locally self-similar Gaussian process, fractional Brownian motion, Hurst exponent estimation, Bahadur representation of sample quantiles.}

\section{Introduction}

Many naturally occuring phenomena can be effectively modelled using self-similar processes. Among the simplest models, one can consider the fractional Brownian motion introduced in the statistics community by \Citet{Mandelbrot68}. Fractional Brownian motion can be defined as the only centered Gaussian process, denoted by $(X(t))_{t\in \RR}$, with stationary increments and with variance function $v(\cdot)$, given by $v(t)=\sigma^2|t|^{2H}$, for all $t\in \RR$. The fractional Brownian motion is an $H$-self-similar process, that is for all $c>0$, $\left( X(ct) \right)_{t \in \mathbb{R}} \stackrel{d}{=} c^H  \left( X(t) \right)_{t \in \mathbb{R}}$ (where $\stackrel{d}{=}$ means equal in finite-dimensional distributions) with autocovariance function behaving like $\mathcal{O}(|k|^{2H-2})$ as $|k|\to +\infty$. So the discretized increments of the fractional Brownian motion (called the fractional Gaussian noise)  constitute a short-range dependent process, when $H<1/2$, and a long-range dependent process, when $H>1/2$. The index $H$ also characterizes the path regularity since the fractal dimension of the fractional Brownian motion is equal to $D=2-H$. According to the context (long-range dependent processes, self-similar processes,\ldots), a very large variety of estimators of the parameter $H$ has been investigated. The reader is referred to \Citet{Beran94a}, \Citet{Coeurjolly00} or \Citet{Bardet03} for an overview of this problem. Among the most often used estimators we have: methods based on the variogram, on the log$-$periodogram {\it e.g.} \Citet{Geweke83} in the context of long-range dependent processes, maximum likelihood estimator (and Whittle estimator) when the model is parametric {\it e.g. } fractional Gaussian noise, methods based on the wavelet decomposition {\it e.g.} \Citet{Flandrin92} or \Citet{Stoev06} and the references therein, and on discrete filtering studied by \Citet{Kent97}, \Citet{Istas97} and \Citet{Coeurjolly01}. %This paper deals with a robust version of one of the two last methods.
We are mainly interested in the last one, which has several similarities with the wavelet decomposition method. Following \Citet{Constantine94}, \Citet{Kent97}, \Citet{Istas97}, in the case when the process is observed at times $i/n$ for $i=1,\ldots,n$, this method is adapted to a larger class than the fractional Brownian motion, namely the class of centered Gaussian processes with stationary increments that are locally self-similar (at zero). A process $(X(t))_{\scriptstyle t\in \mathbb{R}}$ is said to be locally self-similar (at zero) if its variance function, denoted by $v(\cdot)$, satisfies
\begin{equation} \label{lssgp}
v(t) = \Esp (X(t)^2) = \sigma^2 |t|^{2H}\left(1 + r(t) \right), \qquad \mbox{ with } \; r(t) = \oo{1} \mbox{ as } |t|\to 0,
\end{equation}
for some $0<H<1$. An estimator of~$H$ is derived by using the stationarity of the increments and the local behavior of the variance function. When observing the process at regular subdivisions, the stationarity of the increments is crucial since the method based on discrete filtering (and the one based on the wavelet decomposition) essentially uses the fact that the variance of the increments can be estimated by the sample moment of order~2. We do not believe that this framework could be valid for the estimation of the Hurst exponent of Riemann-Liouville's process, {\it e.g.} \Citet{Alos99} which is an $H$-self-similar centered Gaussian process but with increments satisfying only some kind of local stationarity, see Remark~\ref{RL} for more details.

Let us be more specific on the construction of the wavelet decomposition method, see {\it e.g.} \Citet{Flandrin92}: the authors noticed that the variance of the wavelet coefficient at a scale say $j$ behaves like $2^{j(2H-1)}$. An estimator of $H$ is then derived by regressing the logarithm of sample moment of order~2 at each scale against $\log(j)$ for various scales. This procedure exhibits good properties since it is also proved that the more vanishing moments the wavelet has the observations are more decorrelated. And so asymptotic results are quite easy to obtain. However, \Citet{Stoev06} illustrate the fact that this kind of estimator is very sensitive to additive outliers and to non-stationary artefacts. Therefore, they mainly propose to replace at each scale, the sample moment of order~2, by the sample median of the squared coefficients. This procedure, for which the authors assert that no theoretical result is available, is clearly more robust.

The main objective of this paper is to extend the procedure proposed by \Citet{Stoev06} by deriving semi-parametric estimators of the parameter $H$, using discrete filtering methods, for the class of processes defined by (\ref{lssgp}). The procedure is extended in the sense that we consider either convex combinations of sample quantiles or trimmed-means. Moreover, we provide convergence results. The key-ingredient is a Bahadur representation of sample quantiles obtained in a certain dependence framework. Let $\Vect{Y}=(Y(1),\ldots,Y(n))$ be a vector of $n$ i.i.d. random variables with cumulative distribution function $F$, as well denote by $\qt{}{p}$ and $\qtEstSimp{p}$ the quantile respectively the sample quantile of order $p$. By assuming that  $F^{\prime}(\qt{}{p})>0$ and $F^{\prime \prime}(\qt{}{p})$  exists, Bahadur proved that as $n\to+\infty$,
$$
\qtEstSimp{p} - \qt{}{p} = \frac{p-\cdfEstSimp{p}}{f(\qt{}{p}} + r_n,
$$
with $r_n=\Oas{n^{-3/4} \log(n)^{3/4}}$. Using a law of iterated logarithm's type result, Kiefer obtained the exact rate $n^{-3/4}\log \log(n)^{3/4}$. Extensions of the above results to dependent random variables have been pursued in \Citet{Sen72} for $\phi-$mixing variables, in \Citet{Yoshihara95} for strongly mixing variables, and recently in \Citet{Wu05} for short-range and long-range dependent linear processes, following works of \Citet{Hesse90} and \Citet{Ho96}. Our contribution is to provide a Bahadur representation for sample quantiles in another context that is for non-linear functions of Gaussian processes with correlation function decreasing as $k^{-\alpha}L(k)$ for some $\alpha>0$ and some slowly varying function $L(\cdot)$. The bounds for $r_n$ are obtained under the same assumption as those used by \Citet{Bahadur66}.

The paper is organized as follows. In Section~\ref{sec-notation}, we give some basic notations and some background on discrete filtering. In Section~\ref{sec-estimators}, we derive semi-parametric estimators of the parameter $H$, when a single sample path of a process defined by~(\ref{lssgp}) is observed over a discrete grid of the interval $[0,1]$.
% The different estimators are either based on convex combinations of sample quantiles  or trimmed-means.
%, using discrete filtering methods, for the class of processes defined by (\ref{lssgp}).
%Estimators are based on convex combinations of sample quantiles of $\Vect{g(X^a)}$ which is a function of the series $\Vect{X}$ filtered with $\Vect{a}$, that is
%$$\widehat{\xi}_n(\Vect{p},\Vect{c}, \Vect{g(X^a)}) = \sum_{k=1}^K c_k \widehat{\xi}_n(p_k,\Vect{g(X^a)}). $$
%Two functions are considered $g(\cdot)=|\cdot|^\alpha$ for some $\alpha>0$ and $g(\cdot)=\log|\cdot|$, leading to two different estimators of $H$.
Section~\ref{sec-results} presents the main results: Bahadur representations and asymptotic results for our estimators.
% We first establish a Bahadur representation of sample quantiles for non-linear functions of Gaussian sequences and a uniform version of this representation. % sequences with correlation function decreasing hyperbollically.
%These results are then applied to control almost surely our estimators and to obtain, under certain conditions their asymptotic normality with rate $1/\sqrt{n}$, for all $0<H<1$.
In Section~\ref{sec-sim} are presented some numerical computations to compare the theoretical asymptotic variance of our estimators and a simulation study is also given. In particular, we illustrate the relative efficiency with respect to Whittle estimator and the fact that such estimators are more robust than classical ones. Finally, proofs of differents results are presented in Section~\ref{sec-proofs}.

\section{Some notations and some background on discrete filtering} \label{sec-notation}

%Given some standard Gaussian variable with probability (resp. cumulative) distribution function $\phi$ (resp. $\Phi$) denoted $Y$, let us denote for some $0<p<1$, we denote by $\widehat{\xi}_n(p, \Vect{Y})$ (resp. $\mathbb{F}_n(\cdot,\Vect{Y})$) the sample quantile of order $p$ of a vector $\Vect{Y}$ (resp. the sample cumulative distribution function obtained from $\Vect{Y}$).

Given some random variable $Y$, $\cdf{Y}{\cdot}$ denotes the cumulative distribution function of $Y$ and $\qt{Y}{p}$ the quantile of order $p$, $0<p<1$. If \cdf{Y}{\cdot} is absolutely continuous with respect to Lebesgue measure, the probability density function is denoted by \pdf{Y}{\cdot}. The cumulative distribution (resp. probability density) function of a standard Gaussian variable is denoted by $\Phi(\cdot)$ (resp. $\phi(\cdot)$). Based on the observation of a vector $\Vect{Y}=\left(Y(1),\ldots,Y(n) \right)$ of $n$ random variables distributed as $Y$, the sample cumulative distribution function and the sample quantile of order $p$ are respectively denoted by \cdfEst{Y}{\cdot}{Y} and \qtEst{Y}{p}{Y} or simply by \cdfEst{}{\cdot}{Y} and \qtEst{}{p}{Y}. Finally, for some measurable function $g(\cdot)$, we denote by $\Vect{g(Y)}$ the vector of length $n$ with real components $g(Y(i))$, for $i=1,\ldots, n$.

A sequence of real numbers $u_n$ is said to be $\OO{v_n}$ (resp. $\oo{v_n}$) for an other sequence of real numbers $v_n$, if $u_n/v_n$ is bounded (resp. converges to 0 as $n\to +\infty$). A sequence of random variables $U_n$ is said to be $\Oas{v_n}$ (resp. $\oas{v_n}$) if $U_n/v_n$ is almost surely bounded (resp. if $U_n/v_n$ converges towards 0 with probability~1).

%For any random variable $Y$, we also denote by $F(\cdot,Y)$ its cumulative distribution function and by $\xi(p,Y)$ the theoretical quantile of order $p$ associated to $F(\cdot,Y)$. Finally, we denote by $\Vect{g(Y)}$, for some measurable function $g(\cdot)$ the vector of length $n$ with components $g(Y(i))$ for $i=1,\ldots,n$. \\

The statistical model corresponds to a discretized version $\Vect{X}=\left(X(i/n)\right)_{i=1,\ldots,n}$ of a locally self-similar Gaussian process defined by~(\ref{lssgp}).

One of the ideas of our method is to construct some estimators by using some properties of the variance of the increments of $\Vect{X}$ or the variance of the increments of order~2 of $\Vect{X}$. While considering the increments of $\Vect{X}$ is conventional since the associated sequence is stationary, considering the increments of order~2 (or of a higher order) could be stranger. However, the main interest relies upon the fact that the observations of the latter resulting sequences are less correlated than those of the simple increments' sequence. All these vectors can actually be seen as special discrete filtering of the vector $\Vect{X}$. Let us now specify some general background on discrete filtering and its consequence on the correlation structure. The vector $\Vect{a}$ is a filter of length $\ell+1$ and of order $\nu\geq 1$ with real components if
$$ \sum_{q=0}^\ell q^j a_q=0, \mbox{ for } j=0,\ldots,\nu-1 \;\; \quad  \mbox{ and } \quad \sum_{q=0}^\ell q^\nu a_q \neq 0. $$
For example, $\Vect{a}=(1,-1)$ (resp. $\Vect{a}=(1,-2,1)$) is a filter with order 1 (resp. 2). Let $\Vect{X}^{\Vect{a}}$ be the series obtained by filtering $\Vect{X}$ with $\Vect{a}$, then:
$$
X^{\Vect{a}}\left( \frac{i}{n} \right) = \sum_{q=0}^\ell a_q X\left( \frac{i-q}{n} \right) \quad \mbox{ for } i\geq \ell+1.
$$
Applying in turn the filter $\Vect{a}=(1,-1)$ and $\Vect{a}=(1,-2,1)$ leads to the increments of $\Vect{X}$, respectively the increments of $\Vect{X}$ of order~2. One may also consider other filters such as Daubechies wavelet filters, {\it e.g. } \Citet{Daubechies92}.

\noindent The following assumption is needed by different results presented hereafter: \\
$\bm{Assumption \; A_1(k)} :$ for $i=1,\ldots,k$
$$ v^{(i)}(t) = \sigma^2 \beta(i) |t|^{2H-i} + \oo{|t|^{2H-i}} %\Longleftrightarrow r^{(i)}(t)=o(|t|^{2H-i}),
$$
with $\beta(i)=2H(2H-1)\ldots(2H-i+1)$ (where $k\geq 1$ is an integer).
%\item[$\bm{Assumption \; A_2(p)}$:]  $g:\RR\rightarrow \RR$  is a continuous function such that $g(\cdot)=g_1(|\cdot | )$ with $g_1(\cdot)$ bijective on $\RR^+$ and differentiable in a neighbourhood of $\xi(p,g(Y))$.

This assumption assures that the variance function $v(\cdot)$ is sufficiently smooth around~0. It allows us to assert that the correlation structure of a locally self-similar discretized and filtered Gaussian process can be compared to the one of the fractional Brownian motion. This is announced more precisely in the following Lemma.

\begin{lemma}[e.g. \Citet{Kent97}] \label{lemKent} Let $\Vect{a}$ and $\Vect{a^\prime}$ be two filters of length $\ell+1$ and $\ell^\prime+1$, of order $\nu$ and $\nu^\prime \geq 1$. Then we have:
\begin{eqnarray}
\Esp \left( X^{\Vect{a} }\left(\frac{i}{n} \right) X^{\Vect{a^\prime}}\left(\frac{i+j}{n}\right)\right)   & =&
\frac{-\sigma^2}{2} \sum_{q,q^\prime=0}^\ell a_q a_{q^\prime}^\prime  v\left(\frac{q-q^\prime+j}{n}\right)  \nonumber \\
&=& \gamma_n^{\Vect{a},\Vect{a^\prime}}(j) \left( 1 + \delta_n^{\Vect{a},\Vect{a^\prime}}(j) \right), \label{covXa}
\end{eqnarray}
with
\begin{equation} \label{gammaa}
\gamma_n^{\Vect{a},\Vect{a^\prime}}(j) = \frac{\sigma^2}{n^{2H}} \gamma^{\Vect{a},\Vect{a^\prime}}(j) ,  \quad
\gamma^{\Vect{a},\Vect{a^\prime}}(j) =-\frac12\sum_{q,q^\prime=0}^\ell a_q a_{q^\prime}^\prime |q-q^\prime+j|^{2H}
\end{equation}
and
\begin{equation} \label{deltana}
\delta_n^{\Vect{a},\Vect{a^\prime}}(j) = \frac{\sum_{q,q^\prime}a_q a_{q^\prime} |q-q^\prime+j|^{2H} \times r\left(\frac{q-q^\prime+j}{n}\right) }{\gamma^{\Vect{a},\Vect{a^\prime}}(j)}.
\end{equation}
Moreover, as $|j| \to +\infty$
\begin{equation} \label{compGammaa}
\gamma^{\Vect{a},\Vect{a^\prime}}(j) = \OO{ \frac{1}{|j|^{2H-\nu-\nu^\prime}}}.
\end{equation}
Finally, under Assumption $\bm{A_1(\nu+\nu^\prime)}$, as $n \to +\infty$
\begin{equation} \label{compDeltana}
\delta_n^{\Vect{a},\Vect{a^\prime}}(j) = \oo{1}.
\end{equation}
\end{lemma}

\begin{remark}
In the case of the fractional Brownian motion the sequence $\delta_n$ is equal to~0, whereas it converges towards~0 for more general locally self-similar Gaussian processes, such as the Gaussian processes with stationary increments and with variance function $v(t)=1-\exp(-|t|^{2H})$ or $v(t)=\log ( 1+|t|^{2H})$ for which Assumption $\bm{A_1(k)}$ is satisfied (for every $k \geq 1$).
\end{remark}

\begin{remark}\label{RL}
The stationarity of the increments and the local self-similarity required on the process $X(\cdot)$ are important, if the process is observed at times $i/n$ for $i=1,\ldots,n$.  The crucial result of Lemma~\ref{lemKent} is that the variance function of the filtered series behaves asymptotically as $\gamma_n^{\Vect{a}}(0)$. It seems to be difficult to relax the constraint of stationarity. Consider for example the Riemann-Liouville's process, {\it e.g.} \Citet{Alos99}. This process is a Gaussian process which is $H$-self similar Gaussian but with increments satisfying only some kind of local stationarity. Following the computations of \Citet{Lim01}, the variance of the increments' series of the Riemann-Liouville's process is equal to
$$
\Esp\left( \left( X\left(\frac{i+1}n\right) - X\left( \frac in\right) \right)^2\right)=
\frac1{n^{2H}} \frac1{\Gamma(H+1/2)^2} \left\{ I + \frac1{2H}\right\},
$$
with $I=\int_0^i \left( (1+u)^{H-1/2} -u^{H-1/2}\right)^3du + \int_0^{i/n} u^{2H-1}du.$
This integral cannot be asymptotically independent of time. Note that this could be the case if the process is observed at irregular subdivisions. This question has not been investigated.
\end{remark}

Define $\Vect{Y}^{\Vect{a}}$ as the normalized vector $\Vect{X}^{\Vect{a}}$ with variance $1$.
The covariance between $Y^{\Vect{a}}(i/n)$ and  $Y^{\Vect{a^\prime}}(i+j/n)$ is denoted by $\rho_n^{\Vect{a},\Vect{a^\prime}}(j)$.
Under Assumption $\bm{A_1(\nu+\nu^\prime)}$, the following equivalence holds as $n \to +\infty$
\begin{equation} \label{defRhoa}
\rho_n^{\Vect{a},\Vect{a^\prime}}(j) \sim  \rho^{\Vect{a},\Vect{a^\prime}}(j) = \frac{ \gamma^{\Vect{a},\Vect{a^\prime}}(j)}{ \sqrt{\gamma^{\Vect{a},\Vect{a}}(0) \gamma^{\Vect{a^\prime},\Vect{a^\prime}}(0) }}.
\end{equation}
When $\Vect{a}=\Vect{a^\prime}$, we set, for the sake of simplicity $\gamma_n^{\Vect{a}}(\cdot)=\gamma_n^{\Vect{a},\Vect{a}}(\cdot)$, $\delta_n^{\Vect{a}}(\cdot)=\delta_n^{\Vect{a},\Vect{a}}(\cdot)$, $\rho_n^{\Vect{a}}(\cdot)=\rho_n^{\Vect{a},\Vect{a}}(\cdot)$, $\gamma^{\Vect{a},\Vect{a}}(\cdot)=\gamma^{\Vect{a}}(\cdot)$ and $\rho^{\Vect{a}}(\cdot)=\rho^{\Vect{a},\Vect{a}}(\cdot)$.

\section{New estimators of $H$} \label{sec-estimators}

\subsection{Estimators based on a convex combination of sample quantiles}

Let $(\Vect{p},\Vect{c})= (p_k,c_k)_{k=1,\ldots,K} \in ((0,1)\times \RR^+)^K$ for an integer $1\leq K <+\infty$. Define the following statistics based on a convex combination of sample quantiles:
\begin{equation} \label{Lestimate}
\lqtEst{p}{c}{X^a} = \sum_{k=1}^K c_k \; \qtEst{}{p_k}{X^a},
\end{equation}
where $c_k, k=1,\ldots,K$ are positive real numbers such that $\sum_{k=1}^K c_k=1$. For example, this corresponds to the sample median when $K=1, \Vect{p}=1/2, \Vect{c}=1$ , to a mean of quartiles when $K=2, \Vect{p}=(1/4,3/4), \Vect{c}=(1/2,1/2)$ .
% and to the $\beta$ trimmed-mean (for $0<\beta<1$) when $K=[ (1-\beta)n]$ and $p_k=\beta+k/n, c_k=1/ K$ for $k=1,\ldots,K$.
Consider the following computation: from Lemma~\ref{lemKent}, we have, as $n \to +\infty$
$$
\lqtEst{p}{c}{X^a } \sim \frac{\sigma^2}{n^{2H}} \; \gamma^{\Vect{a}}(0) \lqtEst{p}{c}{Y^a }\; .
$$
\begin{remark} \label{unFiltre}
It may be expected that \lqtEst{p}{c}{Y^a } converges towards a constant as $n \to +\infty$. In itself, this result is not interesting, since two parameters remain unknown: $\sigma^2$ and $H$ and thus, it is impossible to derive an estimator of $H$.
\end{remark}
Remark~\ref{unFiltre} suggests that we have to use at least two filters. Among all available filters, let us  consider the sequence $(\Vect{a}^m)_{m\geq 1}$ defined by
$$a_i^m \; = \left\{ \begin{array}{ll}\
a_j&\; \mbox{if } i=jm \\
0  &\; \mbox{otherwise } \\
\end{array} \right. \qquad\qquad \mbox{for } i=0,\ldots,m\ell \;,$$
which is none other than the filter $\Vect{a}$ dilated $m$ times.
For example, if the filter $\Vect{a}=\Vect{a}^1$ corresponds to the filter $(1,-2,1)$, then $\Vect{a}^2=(1,0,-2,0,1)$, $\Vect{a}^3=(1,0,0,-2,0,0,1)$, \ldots As noted by \Citet{Kent97} or \Citet{Istas97}, the filter $\Vect{a}^m$, of length $m\ell+1$, is of order $\nu$ and has the following interesting property~:
\begin{equation} \label{propDil}
\gamma^{\Vect{a}^m}(0)=m^{2H}\gamma^{\Vect{a}}(0).
\end{equation}
From Lemma~\ref{lemKent}, this simply means that $\Esp\left( \Vect{X^{a^m}}(i/n)^2\right)= m^{2H} \Esp\left( \Vect{X^a}(i/n)^2 \right)$, exhibiting some kind of self-similarity property of the filtered coefficients. As specified in the introduction, the same property can be pointed out in the context of wavelet decomposition.

Our methods, that exploit the nice property   (\ref{propDil}), are based on a convex combination of sample quantiles $\lqtEst{p}{c}{g(X^{\Vect{a}^m})}$ for two positive functions $g(\cdot)$: $g(\cdot)=|\cdot|^\alpha$ for $\alpha>0$ and $g(\cdot)=\log|\cdot|$. For such functions $g(\cdot)$ we manage, by using some property established in Lemma~\ref{lemKent}, to define some very simple estimators of the Hurst exponent through a simple linear regression. Other choices of the function $g(\cdot)$ have not been investigated in this paper. At this stage, let us specify that our methods extend the one proposed by \Citet{Stoev06}; indeed they only consider the statistic  $\lqtEst{p}{c}{g(X^{\Vect{a}^m})}$ for $\Vect{p}=1/2$, $\Vect{c}=1$, $g(\cdot)=(\cdot)^2$, that is the sample median of the squared coefficients. From (\ref{gammaa}) and (\ref{propDil}), we have
\begin{eqnarray}
\lqtEst{p}{c}{|X^{\Vect{a}^m}|^\alpha}  &=&
\Esp\left( (X^{\Vect{a^m}}(1/n))^2\right)^{\alpha/2}  \lqtEst{p}{c}{|Y^{\Vect{a}^m}|^\alpha} \nonumber \\
%&=& \frac{\sigma^\alpha}{n^{\alpha H}} \gamma^{\Vect{a}^m}(0)^{\alpha/2}\;
%\widehat{\xi}_n \left( \Vect{p},\Vect{c}, \Vect{|Y^{\Vect{a}^m}|^\alpha} \right) \nonumber \\
&=& m^{\alpha H} \; \frac{\sigma^\alpha}{n^{\alpha H}} \gamma^{\Vect{a}}(0)^{\alpha/2} \; \left( 1 + \delta_n^{a^m}(0)\right)^{\alpha/2}
\lqtEst{p}{c}{|Y^{\Vect{a}^m}|^\alpha} ,\label{idee1}
\end{eqnarray}
and
\begin{eqnarray}
\!\!\!\!\!    \lqtEst{p}{c}{\log|X^{\Vect{a}^m}|} &=& \frac12\log \Esp(X^{\Vect{a^m}}(1/n))^2  +  \lqtEst{p}{c}{\log|Y^{\Vect{a}^m}|}
\nonumber \\
&=& H \log(m) + \log \left( \frac{\sigma^2}{n^{2H}} \gamma^{\Vect{a}}(0)\right) \nonumber \\
&& + \frac{1}{2} \log\left( 1 + \delta_n^{a^m}(0)\right)+ \lqtEst{p}{c}{\log|Y^{\Vect{a}^m}|}\!. \label{idee2}
\end{eqnarray}
Denote by $\kappa_H= n^{-2H} \sigma^2 \gamma^{\Vect{a}}(0)$. Equations~(\ref{idee1}) and~(\ref{idee2}) can be rewritten as
\begin{eqnarray}
\log \lqtEst{p}{c}{|X^{\Vect{a}^m}|^\alpha} &=& \alpha H \log(m) + \log \left( \kappa_H^{\alpha/2} \; \lqt{|Y|^\alpha}{p}{c} \right) + \varepsilon_m^\alpha, \label{idee1Bis}\\
\lqtEst{p}{c}{\log|X^{\Vect{a}^m}|} &=& H \log(m) + \log(\kappa_H )+\; \lqt{\log|Y|}{p}{c} + \varepsilon_m^{\log} \label{idee2Bis}
\end{eqnarray}
with the random variables $\varepsilon_m^\alpha$ and $\varepsilon_m^{\log}$ respectively defined by
\begin{equation} \label{defepsmalpha}
\varepsilon_m^\alpha = \log
\left(  \frac{ \lqtEst{p}{c}{|Y^{\Vect{a}^m}|^\alpha}}{\lqt{|Y|^\alpha}{p}{c}} \right) + \frac{\alpha}{2} \log\left( 1 + \delta_n^{a^m}(0)\right),
\end{equation}
and
\begin{equation} \label{defepsmlog}
\varepsilon_m^{\log}= \lqtEst{p}{c}{\log|Y^{\Vect{a}^m}|} - \lqt{\log|Y|}{p}{c} + \frac{1}{2} \log\left( 1 + \delta_n^{a^m}(0)\right)
\end{equation}
where, for some random variable $Z$, $\lqt{Z}{p}{c}  = \sum_{k=1}^K \; c_k \; \qt{Z}{p_k}$. We decide to rewrite Equations~(\ref{idee1}) and~(\ref{idee2}) as~(\ref{idee1Bis}) and~(\ref{idee2Bis}), since we expect that $\varepsilon_m^\alpha$ and $\varepsilon_m^{\log}$ converge (almost surely) towards~0 as $n \to +\infty$.

From Remark~\ref{unFiltre}, two estimators of $H$ can be defined through a linear regression of $\big(\log \lqtEst{p}{c}{|X^{a^m}|^\alpha}\big)_{m=1,\ldots,M}$
and $\big(\lqtEst{p}{c}{\log|X^{a^m}|}\big)_{m=1,\ldots,M}$
on $\big(\log m\big)_{m=1,\ldots,M}$ for some $M\geq 2$. These estimators are denoted
by $\widehat{H}^\alpha$ and $\widehat{H}^{\log}$.
By denoting $\Vect{A}$ the vector of length $M$ with components $A_m=\log m -\frac1M \sum_{m=1}^M\log(m)$, $m=1,\ldots,M$, we have explicitly from~(\ref{idee1Bis}) and~(\ref{idee2Bis})
and the definition of least squares estimates (see {\it e.g. \Citet{Antoniadis92}}):
\begin{eqnarray}
\widehat{H}^\alpha&=& \frac{\tr{\Vect{A}}}{\alpha||\Vect{A}||^2} \left(  \log
\lqtEst{p}{c}{|X^{\Vect{a}^m}|^\alpha} \right)_{m=1,\ldots,M}, \label{defHalpha} \\
\widehat{H}^{\log} &=& \frac{\tr{\Vect{A}}}{||\Vect{A}||^2} \left(
\lqtEst{p}{c}{\log |X^{\Vect{a^m}}| }     \right)_{m=1,\ldots,M}, \label{defHlog}
\end{eqnarray}
where $||\Vect{z}||$ for some vector $\Vect{z}$ of length $d$ denotes the norm defined by $\left(\sum_{i=1}^d z_i^2 \right)^{1/2}$.

We can point out that $\widehat{H}^\alpha$ and $\widehat{H}^{\log}$ are independent of the scaling coefficient~$\sigma^2$.

\subsection{Estimators based on trimmed means}

Let $0<\beta_1\leq\beta_2<1$ and $\Vect{\beta}=(\beta_1,\beta_2)$, denote by $\tmEst{g(X^a)}$ the $\Vect{\beta}-$trimmed mean of the vector $\Vect{g(X^a)}$ given by
$$
\tmEst{g(X^a)} = \frac1{n-[n\beta_2]-[n\beta_1]}\sum_{[n\beta_1]+1}^{n-[n\beta_2]} \left( \Vect{g(X^a)} \right)_{(i),n},
$$
where $\left( \Vect{g(X^a)} \right)_{(1),n} \leq \left( \Vect{g(X^a)} \right)_{(2),n}\leq \ldots \leq \left( \Vect{g(X^a)} \right)_{(2),n}$ are the order statistics of $\left( \Vect{g(X^a)} \right)_1,\ldots,\left( \Vect{g(X^a)} \right)_n$. It is well-known that $\left(\Vect{g(X^a)}\right)_{(i),n}=\qtEst{}{\frac{i}n}{g(X^a)}$.
Hence, by following the ideas of the previous section, one may obtain
\begin{eqnarray}
\log\left( \tmEst{|X^{a^m}|^\alpha}\right) &=& \alpha H \log(m) + \log \left( \kappa_H^{\alpha/2} \; \tm{|Y|^\alpha} \right) + \varepsilon_m^{\alpha,tm}, \label{idee1BisTM}\\
\tmEst{\log |X^{a^m}|} &=& H \log(m) + \log(\kappa_H) + \tm{\log|Y|} +  \varepsilon_m^{\log,tm} \label{idee2BisTM}
\end{eqnarray}
with
\begin{equation} \label{defepsmalphaTM}
\varepsilon_m^{\alpha,tm}= \tmEst{|Y^{a^m}|^\alpha}-\tm{|Y|^\alpha} +\frac{\alpha}{2} \log\left( 1 + \delta_n^{a^m}(0)\right),
\end{equation}
and
\begin{equation} \label{defepsmlogTM}
\varepsilon_m^{\log,tm}= \tmEst{\log |Y^{a^m}|} -\tm{\log|Y|}+
\frac12\log\left( 1 + \delta_n^{a^m}(0)\right),
\end{equation}
where for some random variable $Z$, $\tm{Z}$ is referring to
\begin{equation} \label{defTMtheorique}
\tm{Z} = \frac1{1-\beta_2-\beta_1} \int_{\beta_1}^{1-\beta_2} \qt{Z}{p} dp.
\end{equation}
%As the previous section, it is hoped that $\varepsilon_m^{\alpha,tm}$ and $\varepsilon_m^{\log,tm}$ converge (almost surely) towards~0. Then, we can define two new estimators of $H$
As in the previous section, two estimators of $H$, denoted by $\widehat{H}^{\alpha,tm}$ and $\widehat{H}^{\log,tm}$, is derived through a log-linear regression
\begin{eqnarray}
\widehat{H}^{\alpha,tm}&=& \frac{\tr{\Vect{A}}}{\alpha||\Vect{A}||^2} \left(   \tmEst{|X^{a^m}|^\alpha} \right)_{m=1,\ldots,M}.\label{defHalphaTM}\\
\widehat{H}^{\log,tm} &=& \frac{\tr{\Vect{A}}}{||\Vect{A}||^2} \left( \tmEst{\log |X^{a^m}|}
\right)_{m=1,\ldots,M}. \label{defHlogTM}
\end{eqnarray}

\begin{remark}
The estimator referred to the ``estimator based on the quadratic variations'' in the simulation study and studied with the same formalism by \Citet{Coeurjolly01} corresponds to the estimator~$\widehat{H}^{\alpha,tm}$ with $\alpha=2$, $\beta_1=\beta_2=0$.
\end{remark}

\section{Main results} \label{sec-results}

To simplify the presentation of different results, consider the two following assumptions on different parameters involved in the estimation procedures

$\bm{Assumption \; A_2(\Vect{p},\Vect{c})} :$  $\Vect{a}$ is a filter of order $\nu \geq 1$, $\alpha$ is a positive  real number, $\Vect{p}$ (resp. $\Vect{c}$) is a vector of length $K$ (for some $1\leq K<+\infty$) such that $0<p_k<1$ (resp. $c_k>0$ and $\sum_{k=1}^K c_k=1$), $M$ is an integer $\geq 2$.

$\bm{Assumption \; A_3(\Vect{\beta})} :$ $\Vect{a}$ is a filter of order $\nu \geq 1$, $\alpha$ is a positive  real number, $\Vect{\beta}=(\beta_1,\beta_2)$ is such that $0<\beta_1\leq \beta_2<1$, $M$ is an integer $\geq 2$.

Since $\tr{A}\left(\log(m)\right)_{m=1,\ldots,M}=||A||^2$ and $\tr{A}\Vect{1}=0$ (where $\Vect{1}=(1)_{m=1,\ldots,M}$), we have
\begin{equation} \label{HestMoinsH}
\widehat{H}^\alpha - H = \frac{\tr{\Vect{A}}}{\alpha ||\Vect{A}||^2} \; \Vect{\varepsilon^\alpha} \quad \mbox{ and } \quad \widehat{H}^{\log} - H = \frac{\tr{\Vect{A}}}{ ||\Vect{A}||^2} \; \Vect{\varepsilon^{\log}},
\end{equation}
and
\begin{equation} \label{HestMoinsHTM}
\widehat{H}^{\alpha,tm} - H = \frac{\tr{\Vect{A}}}{\alpha ||\Vect{A}||^2} \; \Vect{\varepsilon^{\alpha,tm}} \quad \mbox{ and } \quad \widehat{H}^{\log,tm} - H = \frac{\tr{\Vect{A}}}{ ||\Vect{A}||^2} \; \Vect{\varepsilon^{\log,tm}},
\end{equation}
where $\Vect{\varepsilon^\alpha}= \big(  \varepsilon_m^\alpha  \big)_{m=1,\ldots,M}$,
$\Vect{\varepsilon^{\log}}= \big(  \varepsilon_m^{\log}  \big)_{m=1,\ldots,M}$, $\Vect{\varepsilon^{\alpha,tm}}= \big(  \varepsilon_m^{\alpha,tm}  \big)_{m=1,\ldots,M}$ and $\Vect{\varepsilon^{\log,tm}}= \big(  \varepsilon_m^{\log,tm}  \big)_{m=1,\ldots,M}$.
Hence, in order to study the convergence of different estimators, it is sufficient to obtain some convergence results of sample quantiles $\widehat{\xi}(p,\Vect{g(Y^a)})$ for some function $g(\cdot)$ and some filter $\Vect{a}$. Therefore, we first establish a Bahadur representation of sample quantiles for some non-linear function of Gaussian sequences with correlation function decreasing as $k^{-\alpha}$, for some $\alpha>0$. In fact, the existing litterature on nonlinear function of Gaussian sequences ({\it e.g. \Citet{Taqqu77}}) allows us to slighlty extend this framework by considering correlation function decreasing as $k^{-\alpha}L(k)$, for some slowly varying function $L(\cdot)$.

\subsection{Bahadur representation of sample quantiles}

Let us recall some important definitions on Hermite polynomials. The $j$-th Hermite polynomial (for $j\geq 0$) is defined for $t \in \RR$ by
\begin{equation} \label{defHermite}
 H_j(t) = \frac{(-1)^j}{\phi(t)} \; \frac{d^j \phi(t)}{dt^j}.
\end{equation}
%For example, the first polynomials are $H_0(t)=1, H_1(t)=t, H_2(t)=t^2-1, H_3(t)=t^3-3t, \ldots$
The Hermite polynomials form an orthogonal system for the Gaussian measure. More precisely, we have $\Esp\left( H_j(Y) H_k(Y) \right) = j!  \; \delta_{j,k}$. For a measurable function $g(\cdot)$ defined on $\RR$ for which $\Esp(g(Y)^2)<+\infty$, the following expansion holds
$$
g(t) =\sum_{j \geq \tau} \frac{c_j}{j!} \; H_j(t) \quad \mbox{ with } \quad c_j = \Esp\left( g(Y) H_j(Y) \right),
$$
where the integer $\tau$ defined by $\tau= \inf \big\{ j \geq 0, c_j \neq 0 \big\}$,
is called the Hermite rank of the function~$g$. Note that this integer plays an important role. For example, it is related to the correlation of $g(Y_1)$ and $g(Y_2)$ (for $Y_1$ and $Y_2$ two standard gaussian variables with correlation $\rho$) since $\Esp(g(Y_1)g(Y_2))=\sum_{k\geq \tau} \frac{(c_k)^2}{k!} \rho^k\leq \rho^\tau ||g||_{L^2(d\phi)}$.

In order to obtain a Bahadur representation (see {\it e.g.} \Citet{Serfling80}), we have to ensure that $F_{g(Y)}^{\prime}(\qt{}{p})>0$  and $F_{g(Y)}^{\prime \prime}(\cdot)$ exists and is bounded in a neighborhood of $\qt{}{p}$. This is achieved if the function $g(\cdot)$ satisfies the following assumption (see {\it e.g.} \Citet{Dacunha82}, p.33).

$\bm{Assumption \; A_4(\qt{}{p}):}$ there exist $U_i$, $i=1,\ldots,L$, disjoint open sets such that $U_i$ contains a unique solution to the equation $g(t)=\qt{g(Y)}{p}$, such that $F_{g(Y)}^{\prime}(\qt{}{p})>0$ and such that $g$ is a $\mathcal{C}^2-$diffeomorphism on ${\displaystyle \cup_{i=1}^L U_i}$.

\noindent Note that under this assumption
$$
F_{g(Y)}^\prime( \qt{g(Y)}{p} )= \pdf{g(Y)}{\qt{g(Y)}{p}} \; = \;\sum_{i=1}^L \frac{\phi(g_i^{-1}(\qt{}{p}))}{g^\prime(g_i^{-1}(\qt{}{p}))},
$$
where $g_i(\cdot)$ is the restriction of $g(\cdot)$ on $U_i$.
Now, define, for some real $u$, the function $h_u(\cdot)$ by:
\begin{equation} \label{defh}
h_u(t) = \bm{1}_{\{  g(t) \leq u  \}}(t) - \cdf{g(Y)}{u}.
\end{equation}
We denote by $\tau(u)$ the Hermite rank of $h_{u}(\cdot)$. For the sake of simplicity, we set $\tau_p=\tau(\qt{g(Y)}{p})$. For some function $g(\cdot)$ satisfying Assumption $\bm{A_4(\qt{}{p})}$, we denote by
\begin{equation} \label{tauBarre}
\tauB_p = \inf_{\gamma \in \cup_{i=1}^L g(U_i)} \tau(\gamma),
\end{equation}
that is the minimal Hermite rank of  $h_u(\cdot)$ for $u$ in a neighborhood of $\qt{g(Y)}{p}$.

\begin{theorem} \label{bahadur}
Let $\left\{ Y(i) \right\}_{i=1}^{+\infty}$ be a stationary (centered) gaussian process with variance 1, and correlation function $\rho(\cdot)$ such that, as $i\to +\infty$
\begin{equation} \label{hyporho}
| \rho(i) |\sim L(i) \;  i^{-\alpha},
\end{equation}
for some $\alpha>0$ and some slowly varying function at infinity $L(s), s\geq 0$.  Then, under Assumption $\bm{A_4(\qt{}{p})}$, we have almost surely, as $n \to +\infty$
\begin{equation} \label{repBahadur}
\qtEst{}{p}{g(Y)}  - \qt{g(Y)}{p} = \frac{ p-\cdfEst{}{\qt{g(Y)}{p}}{g(Y)}   }{
\pdf{g(Y)}{\qt{g(Y)}{p}}} \; + \; \Oas{r_n(\alpha,\tauB_p) },
\end{equation}
the sequence $\big(r_n(\alpha,\tauB_p)\big)_{n\geq 1}$ being defined by

\begin{equation} \label{defrnAlpha}
r_n (\alpha, \tauB_p) = \left\{ \begin{array}{lll}
n^{-3/4} \log(n)^{3/4} & \mbox{if} & \alpha \tauB_p >1, \\
n^{-3/4} \log(n)^{3/4} L_{\tauB_p}(n)^{3/4} & \mbox{if} & \alpha \tauB_p = 1, \\
n^{-1/2-\alpha\tauB_p/4} \log(n)^{\tauB_p/4+1/2} L(n)^{\tauB_p/4} &\mbox{if}  & 2/3<\alpha \tauB_p<1, \\
n^{-\alpha \tauB_p} \log(n)^{\tauB_p} L(n)^{\tauB_p}& \mbox{if} & 0<\alpha \tauB_p \leq 2/3,
\end{array} \right.
\end{equation}
where for some $\tau\geq 1$, $L_{{\tau}}(n) = \sum_{|i|\leq n}|\rho(i)|^\tau$.
\end{theorem}
Note that if $L(\cdot)$ is an increasing function, $L_\tau(n) =\OO{\log(n) L(n)^\tau}$.

\begin{remark}
Without giving any details here, let us say that the behaviour of the sequence $r_n(\cdot,\cdot)$ is related to the characteristic (short-range or long-range dependence) of the process $\left\{ h_u(Y(i))\right\}_{i=1}^{+^\infty}$ for $u$ in a neighborhood of $\qt{g(Y)}{p}$. In the case $\alpha \tauB_p>1$, corresponding to short-range dependent processes, the result is similar to the one proved by Bahadur, see {\it e.g.} \Citet{Serfling80}, in the i.i.d.  case.
For short-range dependent linear processes, using a law of iterated logarithm's type result \Citet{Wu05} obtained a sharper bound, that is $n^{-3/4} \log \log(n)^{3/4}$. This bound is obtained under the assumption that $F^\prime(\cdot)$ and $F^{\prime \prime}(\cdot)$ exist and are uniformly bounded. For long-range dependent processes ($\alpha \tauB_p \leq 1$), we can observe that the rate of convergence is always lower than $n^{-3/4}\log(n)^{3/4}$ and that the dominant term $n^{-3/4}$ is obtained when $\alpha \tauB_p \to 1$.
\end{remark}

We now propose a uniform Bahadur type representation of sample quantiles. Such a representation has an application in the study of trimmed-mean. For $0<p_0\leq p_1<1$ consider the following assumption which extends $\bm{A_4(\xi(p))}$

$\bm{Assumption \; A_5(p_0,p_1):}$ there exists $U_i$, $i=1,\ldots,L$, disjoint open sets such that $U_i$ contains a solution to the equation $g(t)=\qt{g(Y)}{p}$ for all $p_0\leq p \leq p_1$, such that $F_{g(Y)}^{\prime}(\qt{}{p})>0$ for all $p_0\leq p \leq p_1$ and such that $g$ is a $\mathcal{C}^2-$diffeomorphism on ${\displaystyle \cup_{i=1}^L U_i}$.

\noindent Under the previous assumption, define
\begin{equation} \label{taup0p1}
\tau_{p_0,p_1}= \inf_{\gamma \in \cup_{i=1}^L g(U_i)} \tau(\gamma).
\end{equation}

\begin{theorem} \label{unifBahadur}
Under the conditions of Theorem~\ref{bahadur} and Assumption $\bm{A_5(p_0,p_1)}$, we have almost surely, as $n \to +\infty$
\begin{equation} \label{repUnifBahadur}
\sup_{p_0\leq p \leq p_1} \left| \qtEst{}{p}{g(Y)}  - \qt{g(Y)}{p} - \frac{ p-\cdfEst{}{\qt{g(Y)}{p}}{g(Y)}   }{
\pdf{g(Y)}{\qt{g(Y)}{p}}} \right| = \Oas{r_n(\alpha,\tau_{p_0,p_1})}.
\end{equation}
\end{theorem}

\begin{remark} To obtain convergence results of estimators of $H$, some results are needed concerning sample quantiles of the form $\qtEst{}{p}{g(Y^{a^m})}$, with $g(\cdot)=|\cdot|$. Lemma~\ref{hp} asserts that the Hermite rank $\tau_p$ of the function $h_{\qt{g(Y)}{p}}(\cdot)$ with $g(\cdot)=|\cdot|$, is equal to~2 for all $0<p<1$. Moreover, for all $0<p<1$ and for all $0<p_0\leq p_1 <1$, Assumptions $\bm{A_4(\xi(p))}$ and $\bm{A_5(p_0,p_1)}$ are satisfied, and we have $\tauB_p=\tau_{p_0,p_1}=2$. Since from Lemma~\ref{lemKent}, the correlation function of $\Vect{Y}^{\Vect{a}^m}$ satisfies~(\ref{hyporho}) with $\alpha=2\nu-2H$ and $L(\cdot)=1$, by applying Theorem~\ref{bahadur}, the sequence $r_n(\cdot,\cdot)$ is then given by
\begin{equation} \label{defrnAlphalssgpnu2}
r_n(2\nu-2H, 2)=n^{-3/4} \log(n)^{3/4},\quad  \mbox{ if }\;  \nu\geq 2
\end{equation}
and for $\nu=1$
\begin{equation} \label{defrnAlphalssgpnu1}
r_n (2-2H, 2) = \left\{ \begin{array}{ll}
n^{-3/4} \log(n)^{3/4} & \mbox{if } 0<H< 3/4 ,\\
n^{-3/4} \log(n)^{3/2}  & \mbox{if } H=3/4 , \\
n^{-1/2-(1-H)} \log(n)  &\mbox{if } 3/4<H<5/6, \\
n^{-2(2-2H)} \log(n)^2& \mbox{if }  5/6 \leq H<1 .\\
\end{array} \right.
\end{equation}
\end{remark}

\subsection{Convergence results of estimators of $H$}

In order to specify convergence results, we make the following assumption concerning the remainder term of the variance function $v(\cdot)$.

\noindent $\bm{Assumption \; A_6(\eta)} :$ there exists $\eta>0$ such that $v(t)=\sigma^2|t|^{2H}\left( 1+ \OO{|t|^\eta}\right),$ as $|t|\to 0.$

The first result concentrates itself on estimators $\widehat{H}^{\alpha}$ and $\widehat{H}^{\log}$ based on a convex combination of sample quantiles.\begin{theorem} \label{convHest} Under Assumptions $\bm{A_1(2\nu)}$, $\bm{A_2(p,c)}$ and $\bm{A_6(\eta)}$,

$(i)$ we have almost surely, as $n \to +\infty$
\begin{equation} \label{defyn}
 \widehat{H}^\alpha -H  =
\left\{
\begin{array}{ll}
\OO{n^{-\eta}} + \Oas{n^{-1/2}\log(n)} & \mbox{if } \nu >H+\frac1{4}, \\
\OO{n^{-\eta}} +\Oas{n^{-1/2}\log(n)^{3/2}}& \mbox{if } \nu=1, H= \frac3{4}, \\
\OO{n^{-\eta}} +\Oas{n^{-2(1-H)} \log(n)} & \mbox{if } \nu=1, \frac3{4} <H<1.
\end{array}
\right.
\end{equation}
A similar result holds for $\widehat{H}^{\log}$.

$(ii)$ the mean squared errors (MSE) of $\widehat{H}^{\alpha}$ 
satisfies
\begin{equation} \label{mseHalpha}
\mbox{MSE}\left( \widehat{H}^{\alpha}-H\right)=
\OO{v_n(2\nu-2H)} +
\OO{r_n(2\nu-2H,2)^2} + \OO{n^{-2\eta}}.
\end{equation}
The sequence $r_n(2\nu-2H,2)$ is given by~(\ref{defrnAlphalssgpnu2}) and~(\ref{defrnAlphalssgpnu1}) and the sequence $v_n(\cdot)$ is defined by
\begin{equation} \label{defvn}
v_n(2\nu-2H) = \left\{
\begin{array}{ll}
n^{-1} & \mbox{if } \nu > H+\frac1{4}, \\
n^{-1}\log(n) & \mbox{if } \nu=1, H= \frac34, \\
n^{-4(1-H)} & \mbox{if } \nu=1, \frac34 <H<1.
\end{array} \right.
\end{equation}
Again, the same result holds for $MSE\left( \widehat{H}^{\log}-H\right)$.

$(iii)$ if the filter $\Vect{a}$ is such that $\nu>H+1/4$, and if $\eta>1/2$,
 then we have the following convergence in distribution, as $n \to +\infty$
\begin{equation} \label{convHEstLoi}
\sqrt{n} \left( \widehat{H}^{\alpha}  - H  \right) \longrightarrow \mathcal{N}(0, \sigma^2_\alpha ) \qquad \mbox{ and } \qquad \sqrt{n} \left( \widehat{H}^{\log}  - H  \right) \longrightarrow \mathcal{N}(0, \sigma^2_0 ),
\end{equation}
where $\sigma^2_\alpha$ is defined for $\alpha\geq 0$ by
\begin{equation} \label{sig2Alpha}
\sigma^2_\alpha = \sum_{i \in \ZZ} \sum_{j \geq 1} \frac{1}{(2j)!}  \bigg(
\sum_{k=1}^K \frac{H_{2j-1}(q_k) c_k}{q_k} \; {\pi_k^\alpha}
\bigg)^2 \tr{\Vect{B}} \Mat{R}(i,j) \Vect{B}.
\end{equation}
The vector $\Vect{B}$ is defined by ${\displaystyle \Vect{B}=\frac{\tr{\Vect{A}}}{||\Vect{A}||^2}}$, and the real numbers $q_k$ and  $\pi_k^\alpha$ are defined  by
\begin{equation} \label{defPikAlpha}
q_k=\Phi^{-1} \left(\frac{1+p_k}{2}\right) \quad \mbox{ and } \quad \pi_k^\alpha = \frac{  (q_k)^\alpha    }{ \sum_{j=1}^K c_j  (q_j)^\alpha }.
\end{equation}
Finally, the matrix $\Mat{R}(i,j)$, defined for $i \in \ZZ$ and $j \geq 1$, is a $M \times M$ matrix whose $(m_1,m_2)$ entry is
\begin{equation}\label{defMatR}
\left( \Mat{R}(i,j) \right)_{m_1,m_2} = \rho^{\Vect{a^{m_1}}, \Vect{a^{m_2}} }(i)^{2j},
\end{equation}
where $\rho^{\Vect{a^{m_1}}, \Vect{a^{m_2}} }(\cdot)$ is the correlation function defined by (\ref{defRhoa}).
\end{theorem}
\bigskip

\begin{remark}
The expression of the variance $\sigma_\alpha^2$ given by (\ref{sig2Alpha}) could appear to be very complicated. However, given some vectors $\Vect{p}$ and $\Vect{c}$ and some integer $M$, it does not take unreasonnable effort to compute it for each value of $H$ by truncating the two series. This issue is investigated in Section~\ref{sec-sim} to compare the different parameters.\end{remark}

\begin{remark} \label{rem-mse}
Let us discuss the result (\ref{mseHalpha}). The first term, $\OO{v_n}$, is due to the variance of the sample cumulative distribution function. The second term, $\OO{r_n^2}$ is due to the departure of $\qtEstSimp{p}-\qt{}{p}$ from $\cdfEstSimp{\qt{}{p}}-p$. We leave the reader to check that
$$
\OO{r_n(2\nu -2H,2)^2}  + \OO{v_n(2\nu-2H)} = \left\{
\begin{array}{ll}
\OO{v_n(2\nu-2H)} & \mbox{ if } \nu \geq H+\frac1{4}, \\
\OO{r_n(2\nu -2H,2)^2} &\mbox{ if } \nu < H+\frac1{4}.
\end{array} \right.
$$
Finally, the third one, $\OO{n^{{-2\eta}}}$ is a bias term due to the misspecification of the variance function $v(\cdot)$ around 0.
\end{remark}

\begin{remark} \label{rem-normalite}
If $K=1$, we have, for every $\alpha>0$,
$$\sigma^2_\alpha=\sigma^2_0 = \sum_{i \in \ZZ} \; \sum_{j \geq 1}  \frac{H_{2j-1}(q)^2 }{q^2 \; (2j)!} \;  \tr{\Vect{B}} \Mat{R}(i,j) \Vect{B}.
$$
Assume $\bm{A_6(\eta)}$ with $\eta>1/2$ which allows to neglict the bias term with respect to the variance one. The result (\ref{convHEstLoi}) is proved by using some general central limit theorem obtained in this dependence context by \Citet{Arcones94}, which is available as soon as $\rho^{\Vect{a}}(\cdot)^2$ is summable. Therefore, if only $\bm{A_1(2)}$ is assumed, the filter $\Vect{a}$ cannot exceed~1 (and then correspond to $\Vect{a}=(1,-1)$) and, due to~(\ref{compGammaa}), the result (\ref{convHEstLoi}) is valid only for $0<H<3/4$. As a practical point of view, one observes that for such a filter and large values of $H$, the estimators have very big variance. Note that if $\bm{A_1(2\nu)}$ can be assumed for $\nu>1$, then the asymptotic normality is valid for all the values of $H$.
\end{remark}

\bigskip

\noindent The next result asserts the link between $\widehat{H}^{\log}$ and $\widehat{H}^{\alpha}$.

\begin{corollary} \label{corolHest} Let $(\alpha_n)_{n\geq 1}$ be a sequence such that $\alpha_n \to 0$, as $n \to +\infty$. Then,  under conditions of Theorem~\ref{convHest}~$(ii)$, the following convergence in distribution holds, as $n \to +\infty$
\begin{equation} \label{convHAlphan}
\sqrt{n} \left( \widehat{H}_n^{\alpha_n}-H\right) \longrightarrow \mathcal{N}(0,\sigma^2_0).
\end{equation}
\end{corollary}

The following theorem presents the analog results obtained for the estimators $\widehat{H}^{\alpha,tm}$ and $\widehat{H}^{\log,tm}$ based on trimmed-means.
\begin{theorem} \label{convHestTM}
Under Assumptions  $\bm{A_1(2\nu)}$, $\bm{A_3(\Vect{\beta})}$ and $\bm{A_6(\eta)}$, properties $(i)$ and $(ii)$ of Theorem~\ref{convHest} hold for the estimator $\widehat{H}^{\alpha,tm}$ and $\widehat{H}^{\log,tm}$ with the same rates of convergences. 

$(iii)$ if the filter $\Vect{a}$ is such that $\nu>H+1/4$ and if $\eta>1/2$,
 then, under the notations of Theorem~\ref{convHest}, we have the following convergence in distribution, as $n \to +\infty$
\begin{equation} \label{convHEstLoiTM}
\sqrt{n} \left( \widehat{H}^{\alpha,tm}  - H  \right) \longrightarrow \mathcal{N}(0, \sigma^2_{\alpha,tm} ) \quad \mbox{ and } \quad \sqrt{n} \left( \widehat{H}^{\log,tm}  - H  \right) \longrightarrow \mathcal{N}(0, \sigma^2_{0,tm} ),
\end{equation}
where $\sigma^2_{\alpha,tm}$ is defined for $\alpha\geq 0$ by
\begin{equation} \label{sigma2TM}
\sigma^2_{\alpha,tm}= \sum_{i \in \ZZ} \sum_{j\geq 1} \frac1{(2j)!} \left(
\frac{ \int_{\beta_1}^{1-\beta_2} H_{2j-1}(q) q^{\alpha-1} dp}{\int_{\beta_1}^{1-\beta_2} q^\alpha dp}
\right)^2 \tr{\Vect{B}} \Mat{R}(i,j) \Vect{B},
\end{equation}
with $q=\Phi^{-1}\left(\frac{1+p}2 \right)$.
\end{theorem}

\section{Numerical computation and simulations} \label{sec-sim}

\subsection{Asymptotic constants $\sigma^2_{\alpha}$ and $\sigma^2_{\alpha,tm}$}

\begin{figure}[H]
\hspace*{-1cm}\vspace*{-.2cm}\begin{tabular}{ll}
\includegraphics[scale=.36]{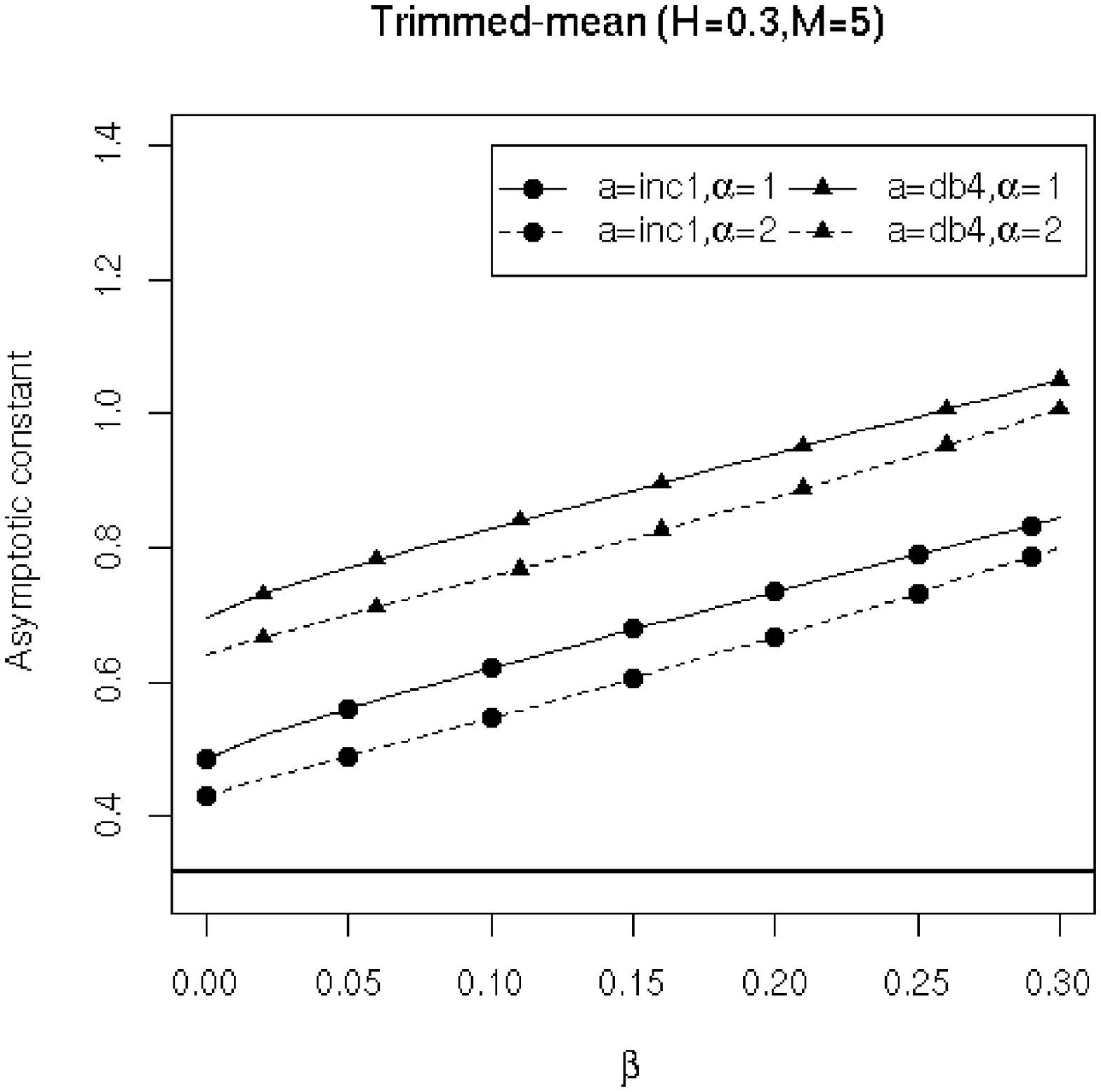} &\includegraphics[scale=.36]{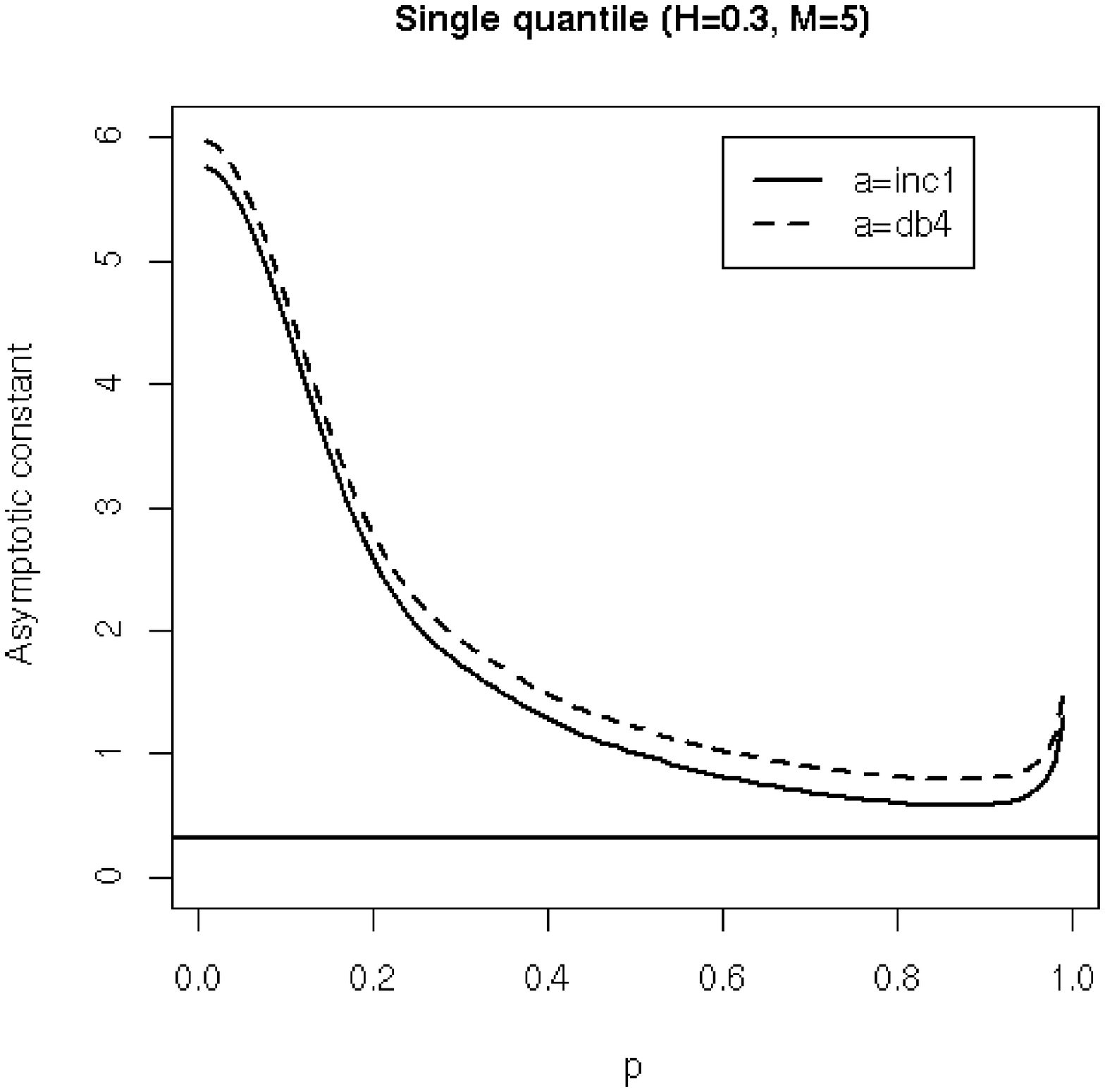} \\
\includegraphics[scale=.36]{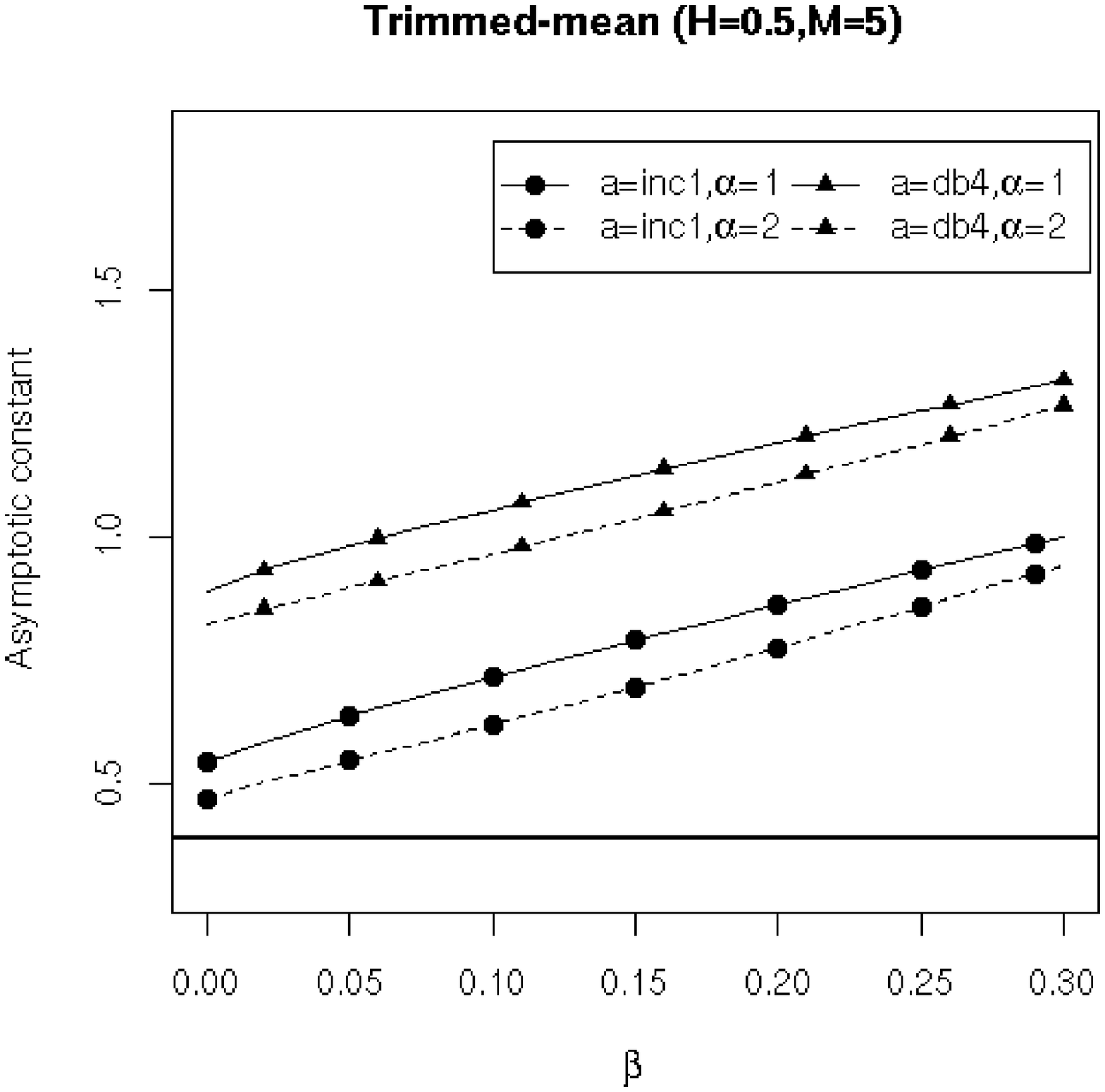} &\includegraphics[scale=.36]{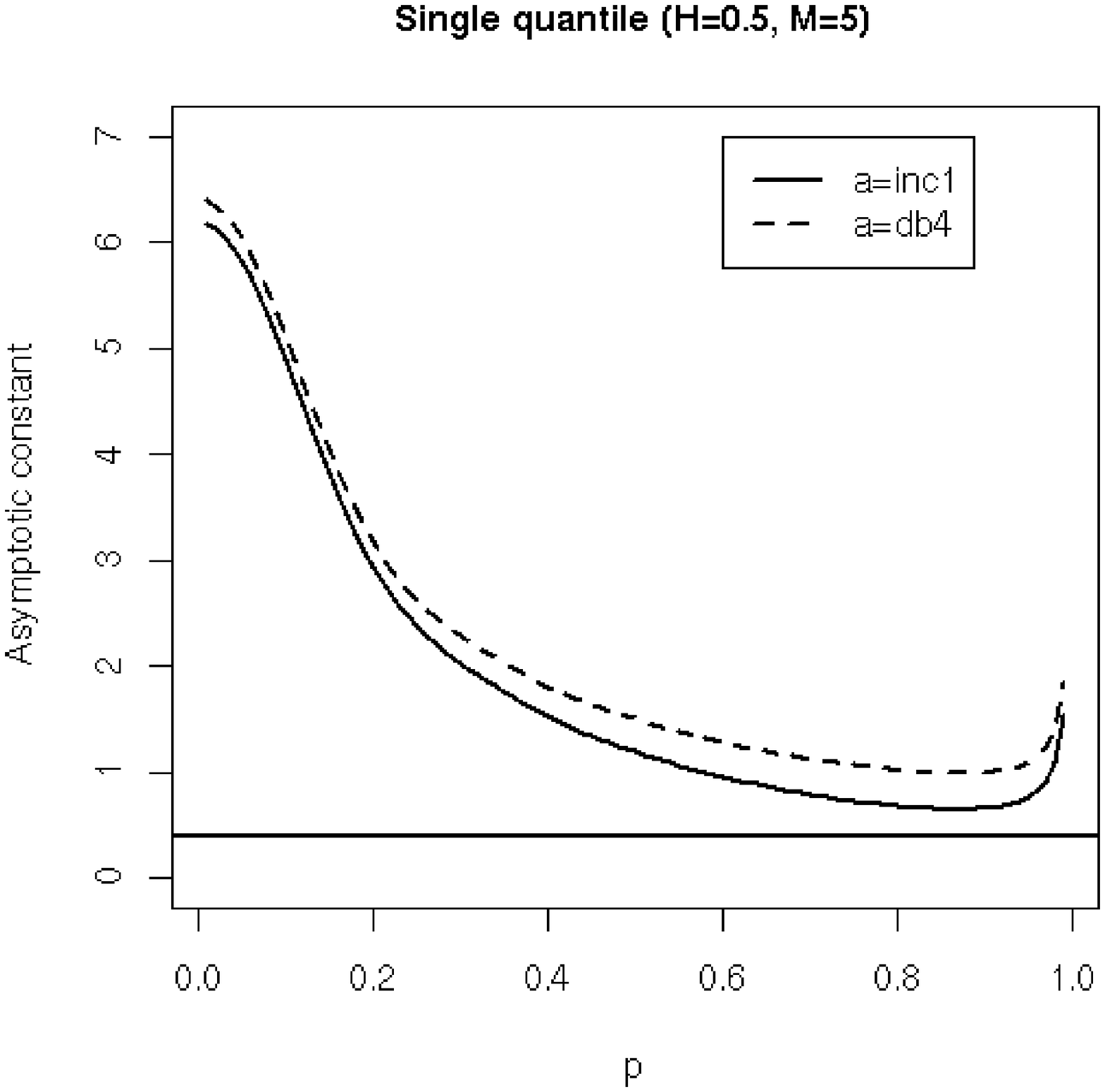}\\
\includegraphics[scale=.36]{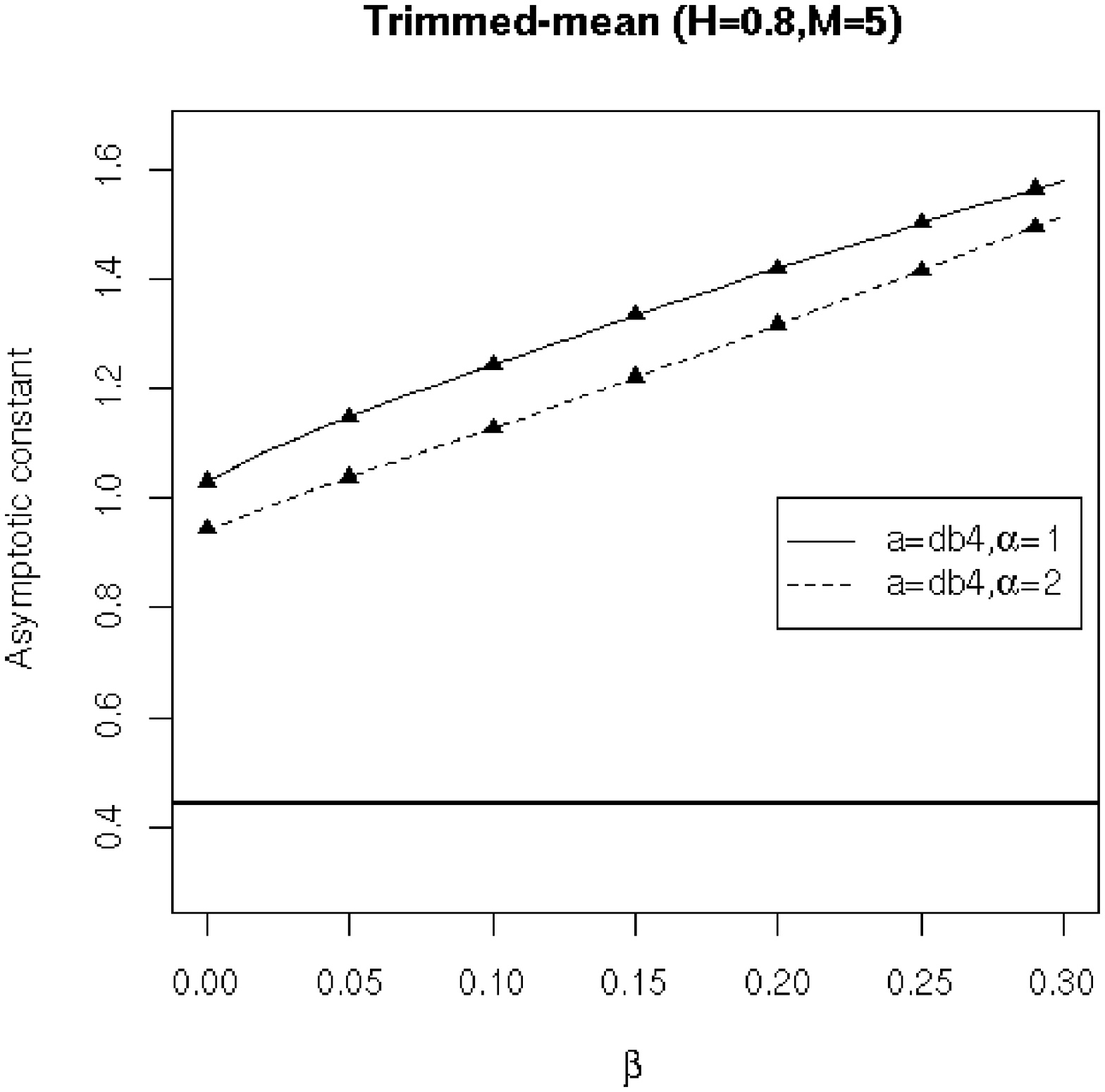} &\includegraphics[scale=.36]{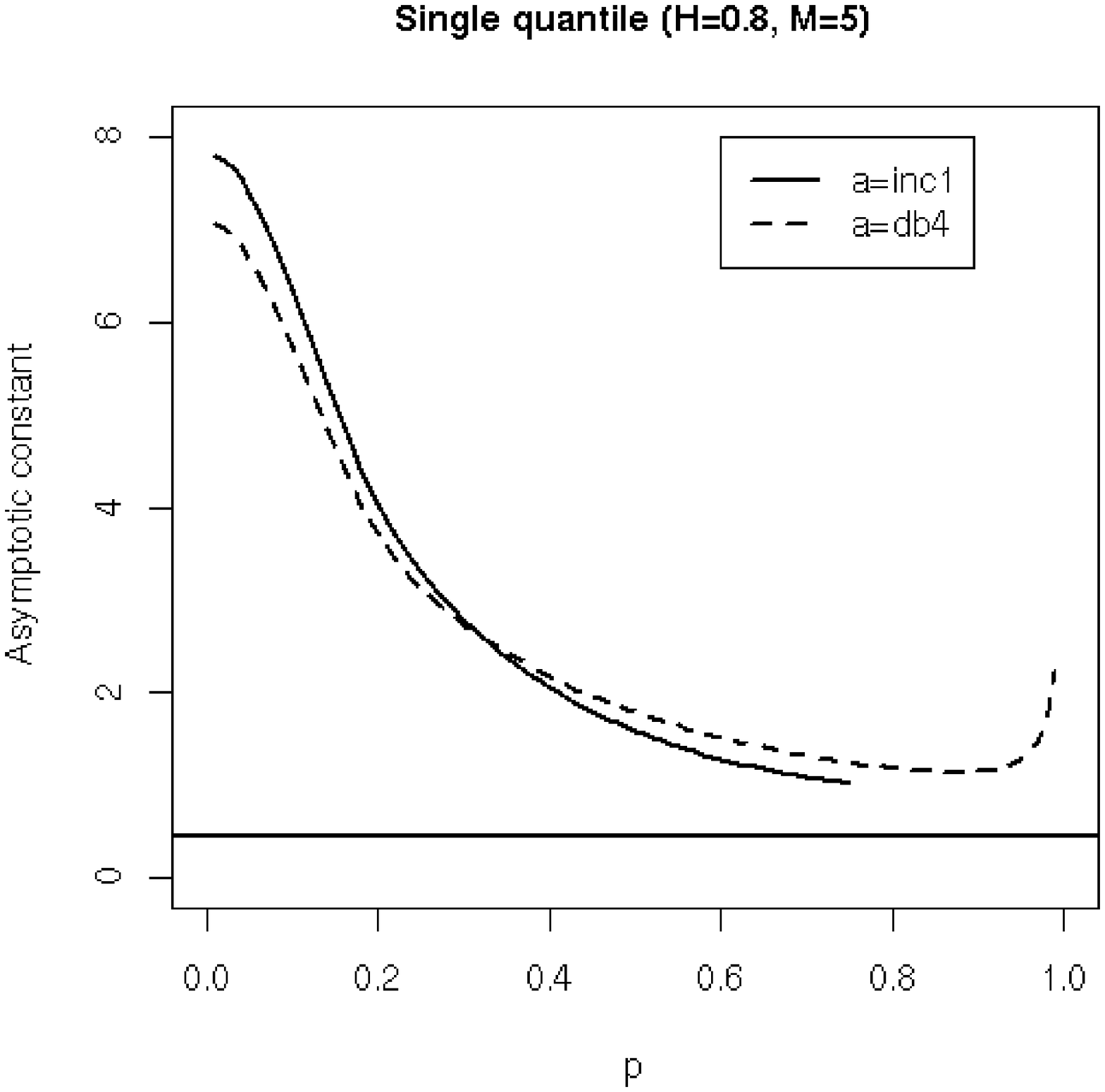}\\
\end{tabular}
\caption{Left: $\sigma^2_{\alpha,tm}$ in terms of $\Vect{\beta}$;  Right: $\sigma^2_{\alpha}$ for estimators based on a single quantile in terms of $p$. Three values of the parameter $H$ are considered: 0.3 (top), 0.5 (middle), 0.8 (bottom). The parameter $M$ is fixed to $M=5$. The constant line corresponds to the asymptotic variance of the Whittle's estimator.}
\label{csteSingleTM}
\end{figure}

In order to compare the different estimators, we intend to compute the asymptotic constants $\sigma^2_{\alpha}$ and $\sigma^2_{\alpha,tm}$ defined by~(\ref{sig2Alpha}) and~(\ref{sigma2TM}) for various set of parameters ($\Vect{a},\Vect{p},\Vect{c}, \Vect{\beta},M$). For this work, both series defining $\sigma^2_{\alpha}$ and $\sigma^2_{\alpha,tm}$ are truncated ($|i|\leq200$, $j\leq 150$).
Figure~\ref{csteSingleTM} illustrates a part of this work. We can propose the following general remarks:

$\bullet$ Among all filters tested, the best one seems to be
$$
\Vect{a}^\star = \left\{
\begin{array}{ll}
inc1 & \mbox{ if } 0<H<3/4,\\
%inc1 & \mbox{ if } H<7/8 \mbox{ and } K=1, p=2\Phi(\sqrt{3})-1, \\
db4 & \mbox{ otherwise.}
\end{array} \right.
$$
where $inc1$ and $db4$ respectively denote the filter $(1,-1)$ and the Daubechies wavelet filter with two zero moments explicitly given by
$$
db4=\left( 0.4829629,-0.8365763,0.22414386,0.12940952 \right).
$$

$\bullet$ Choice of $M$: increasing $M$ seems to reduce the asymptotic constant $\sigma_\alpha^2$. Obviously, a too large $M$ increases the bias since
$\lqtEst{p}{c}{g(X^{a^M})}$ or $\tmEst{g(X^{a^M})}$ are estimated with $N-M\ell$ observations. We recommend setting it to the value~5.

$\bullet$
We did not manage (theoretically and numerically since series defining (\ref{sig2Alpha}) and~(\ref{sigma2TM}) are truncated) to determine the optimal value of $\alpha$. However, for examples considered, it should be near the value $2$.

$\bullet$ Again, this is quite difficult to know theoretically and numerically which choice of $\Vect{p}$ is optimal. What we observed  is that, for fixed parameters $\Vect{a}$, $M$ and $\alpha$, the asymptotic constants are very close to each other.

$\bullet$ Choice of $p$ in the case of a single quantile (see Figure~\ref{csteSingleTM}): the optimal $p$ seems to be near the value $90\%$. However, $p=1/2$, corresponding to the estimator based on the median, leads to good results.

$\bullet$ Choice of $\beta_1=\beta_2=\beta$ for the estimators based on trimmed-means (see Figure~\ref{csteSingleTM}): obviously the constant grows with $\beta$ but we can point out that estimators based on $10\%-$trimmed-means are very competitive with the ones obtained by quadratic variations ($\beta=0$).

\subsection{Simulation}

\begin{figure}[H]
\begin{center}
\includegraphics[scale=.4]{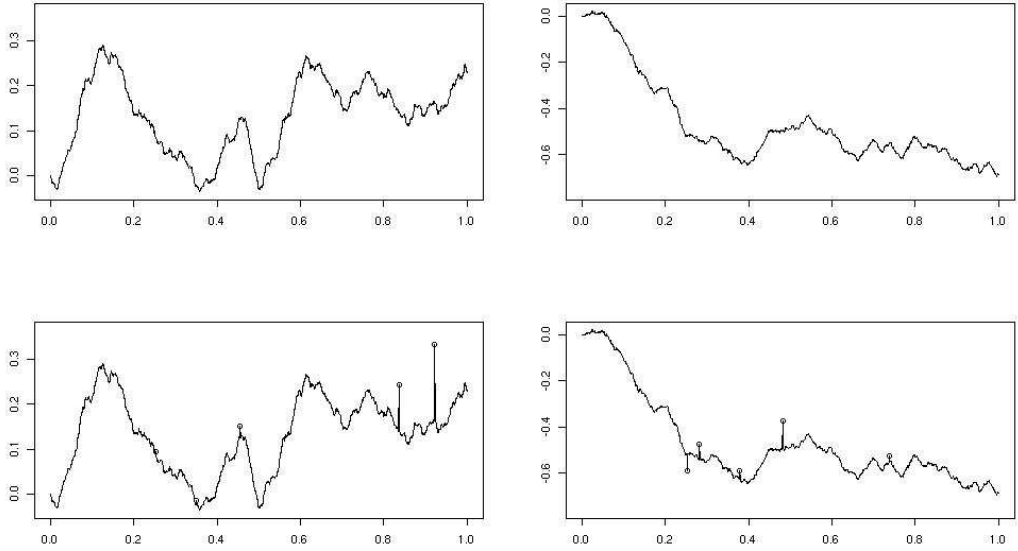} \\
\end{center}
\caption{\small Two examples for the sample paths of non-contaminted (top) and
contaminated processes
with variance function $v(\cdot)=|\cdot|^{2H}$ (left), respectively $v(\cdot)=1-\exp(-|\cdot|^{2H})$ (right), see~(\ref{contamination}).}
\label{exRob}
\end{figure}

A short simulation study is proposed in Table~\ref{tab-Robuste} and Figure~\ref{exRob} for $n=1000$ and $H=0.8$. We consider two locally self-similar Gaussian processes whose variance functions are in turn $v(t)=|t|^{2H}$ (fractional Brownian motion) and $v(t)=1-\exp(-|t|^{2H})$. To generate sample paths discretized over a grid $[0,1]$, we use the method of circulant matrix (see \Citet{Wood94}), which is particularly fast, even for large sample sizes. Various versions of estimators are considered and compared with classical ones, that is the one based on quadratic variations, \Citet{Coeurjolly01}, and the Whittle estimator, \Citet{Beran94a}. In order to illustrate the robustness of our estimators, we also applied them to contaminated version of sample path processes. We obtain a new sample path
 discretized at times $i/n$ and denoted by $X^C(i/n)$ for $i=1,\ldots,n$ through the following model
\begin{equation} \label{contamination}
X^C(i/n)= X(i/n) + U(i) V(i),
\end{equation}
where $U(i)$, $i=1,\ldots,n$ are Bernoulli independent variables $\mathcal{B}(0.005)$, and $V(i)$, $i=1,\ldots,n$ are independent centered Gaussian variables with variance $\sigma_C^2(i)$ such that the signal noise ratio at time $i/n$ is equal to 20 dB.
As a general conclusion of Table~\ref{tab-Robuste}, one can say that all versions of our estimators are very competitive with classical ones when the processes are observed without contamination and they seem to be particularly robust to additive outliers. Both bias and variance are approximately unchanged. This is clearly not the case for classical estimators. Indeed, concerning quadratic variations' method, the estimation procedure is based on the estimation of $\Esp( (X^{\Vect{a}^m}(1/n))^2)$ by sample mean of order 2 of $(\Vect{X}^{\Vect{a}^m})^2$, \Citet{Coeurjolly01}), that is particularly sensitive to additive outliers. Bad results of Whittle estimator can be explained by the fact that maximum likelihood methods are also non-robust methods.

\section{Proofs} \label{sec-proofs}

We denote by $||\cdot||_{L^2(d\phi)}$ (resp. $||\cdot||_{\ell^q}$) the norm defined by $||h||_{L^2(d\phi)}=E(h(Y)^2)^{1/2}$ for some measurable function $h(\cdot)$ (resp. $(\sum_{i \in \ZZ} |u_i|^q)^{1/2}$ for some sequence $(u_i)_{i \in \ZZ}$). In order to simplify the presentation of proofs, we use the notations $\cdf{}{\cdot}$, $\qt{}{\cdot}$, $\pdf{}{\cdot}$, $\cdfEstSimp{\cdot}$ and $\qtEstSimp{\cdot}$ instead of $\cdf{g(Y)}{\cdot}$, $\qt{g(Y)}{\cdot}$, $\pdf{g(Y)}{\cdot}$, $\cdfEst{g(Y)}{\cdot}{g(Y)}$ and $\qtEst{g(Y)}{\cdot}{g(Y)}$ respectively. For some real $x$, $[x]$ denotes the integer part of $x$. Finally, $\lambda$ denotes a generic positive constant.

\subsection{Sketch of the proof of Theorem~\ref{bahadur}} \label{intro-proof}

We give here a brief explanation of the strategy to prove~Theorem~\ref{bahadur}. This proof follows exactly the one proposed by \Citet{Serfling80} in the i.i.d. case. One starts by writing
$$
\frac{p- \cdfEstSimp{\qt{}{p}}}{f(\qt{}{p})} -\left( \qtEstSimp{p}-\qt{}{p}\right)  = A(p) + B(p) + C(p),$$
with
\begin{eqnarray}
A(p)&=&\frac{ p-\cdfEstSimp{\qtEstSimp{p}}}{f(\qt{}{p})} \label{defA}\\
B(p)&=& \frac{ \cdfEstSimp{\qtEstSimp{p}} - \cdfEstSimp{\qt{}{p}} -\left( \cdf{}{\qtEstSimp{p}} -\cdf{}{\qt{}{p}} \right)}{f(\qt{}{p})} \label{defB}\\
C(p)&=& \frac{\cdf{}{\qtEstSimp{p}} -\cdf{}{\qt{}{p}} }{f(\qt{}{p})} - \left(\qtEstSimp{p}-\qt{}{p}\right). \label{defC}
\end{eqnarray}
From the definition of sample quantile, we have almost surely, see {\it e.g.} \Citet{Serfling80}, $A(p)= \Oas{n^{-1}}$. Now, in order to control the term $C(p)$, Taylor's Theorem is used and a control of $\qtEstSimp{p}-\qt{}{p}$ is needed. The latter one is done by Lemma~\ref{lemme-controleQt} which exhibits the sequence $\varepsilon_n(\alpha,\tau_p)$ such that $\qtEstSimp{p}-\qt{}{p}=\Oas{\varepsilon_n(\alpha,\tau_p)}$. Then, in order to control $B(p)$ it is sufficient to control the random variable
$$
S_n(\qt{}{p},\varepsilon_n(\alpha,\tau_p)) = \sup_{|x|\leq \varepsilon_n(\alpha,\tau_p)}
\left|  \Delta(\qt{}{p}+x) -  \Delta(\qt{}{p})\right|,
$$
with $\Delta(\cdot)=\cdfEstSimp{\cdot}-\cdf{}{\cdot}$. This result is detailed in Lemma~\ref{lemme-bahadur}. In order to specify the rate explicited by Theorem~\ref{bahadur}, we present and prove Lemmas~\ref{lemme-controleQt} and~\ref{lemme-bahadur}. Some preliminary results, given by Lemma~\ref{moyps}, Corollary~\ref{corIneg} and Lemma~\ref{lemme-coeff}, are needed. Among other things, Lemma~\ref{moyps} and Corollary~\ref{corIneg} propose some inequalities for controlling the sample mean of non-linear function of Gaussian sequences with correlation function satisfying~(\ref{hyporho}).

\subsection{Auxiliary Lemmas for the proof of Theorem~\ref{bahadur}}

\begin{lemma} \label{moyps}
Let $\left\{ Y(i) \right\}_{i=1}^{+\infty}$ a gaussian stationary process with variance 1 and correlation function $\rho(\cdot)$ such that, as $i\to +\infty$, $| \rho(i) |\sim L(i) i^{-\alpha}$, for some $\alpha>0$ and some slowly varying function at infinity $L(\cdot)$. Let $h(\cdot) \in L^2\left(d\phi\right)$  and denote by $\tau$ its Hermite rank. Define
$$
\overline{Y}_n = \frac1n \sum_{i=1}^n h(Y(i)).
$$
Then, for all $\gamma>0$, there exists a positive constant $\kappa_\gamma=\kappa_\gamma(\alpha,\tau)$, such that
\begin{equation} \label{moypc}
\Prob{ | \overline{Y}_n | \geq \kappa_\gamma y_n} = \OO{n^{-\gamma}},
\end{equation}
with
\begin{equation} \label{defeps}
y_n =y_n(\alpha,\tau) = \left\{ \begin{array}{lll}
n^{-1/2} \log(n)^{1/2}  & \mbox{if } & \alpha \tau >1 ,\\
n^{-1/2} \log(n)^{1/2} {L_\tau(n)}^{1/2} & \mbox{if } & \alpha \tau =1 ,\\
n^{-\alpha\tau/2} \log(n)^{\tau/2} L(n)^{\tau/2} & \mbox{if } & 0<\alpha \tau <1. \\
\end{array}
\right.
\end{equation}
where $L_\tau(n) = \sum_{|i|\leq n}|\rho(i)|^\tau$. In the case $\alpha\tau=1$, we assume that for all $j>\tau$, the limit,  $\lim_{n\to+\infty} L_\tau(n)^{-1} \sum_{|i|\leq n} |\rho(i)|^j$ exists.
\end{lemma}

\begin{proof}
Let $(y_n)_{n\geq 1}$ be the sequence defined by (\ref{defeps}). The proof is splitted into three parts according to the value of $\alpha \tau$.

%% alpha tau<1
\noindent \underline{\bf Case  $\bm{ \alpha \tau <1}$ :}
From Chebyshev's inequality, we have for all $q \geq 1$
$$\Prob{| \overline{Y}_n | \geq \kappa_\gamma y_n} \leq \frac{1}{\kappa_\gamma^{2q} y_n^{2q}} \Esp \bigg( \big(\overline{Y}_n \big)^{2q} \bigg).
$$
From Theorem 1 of \Citet{Breuer83} and in particular Equation (2.6), we have, as $n \to +\infty$
\begin{equation} \label{contMom}
\Esp \bigg( \big(\overline{Y}_n \big)^{2q} \bigg) \sim \frac{(2q)!}{2^q q!} \frac1{n^q} \sigma^{2q},
\quad \mbox{ with } \quad \sigma^2 = \sum_{i \in \ZZ} \sum_{j \geq \tau} \frac{(c_j)^2}{j!} \rho(i)^{j},
\end{equation}
where $c_j$ denotes the $j$-th Hermite coefficient of $h(\cdot)$. Note that $\sigma^2 \leq ||h||^2_{L^2(d\phi)} ||\rho||^2_{\ell^\tau}$. Thus, for $n$ large enough, we have
\begin{equation} \label{cas1eq1}
\Prob{| \overline{Y}_n | \geq \kappa_\gamma y_n} \leq \frac{\lambda}{n^q y_n^{2q}}  \frac{(2q)!}{2^q q!} 
\left( ||h||^2_{L^2(d\phi)} ||\rho||^2_{\ell^\tau}  \kappa_\gamma^{-2} \right)^q .
\end{equation}
From Stirling's formula, we have as $q \to +\infty$
\begin{equation} \label{eqStirling}
\frac{(2q)!}{2^q q!}\; \sim  \; \sqrt{2}\;  q^q  \; (2e^{-1})^q.
\end{equation}
From (\ref{defeps}) by choosing $q=[\log(n)]$, (\ref{cas1eq1}) becomes
$$
\Prob{| \overline{Y}_n | \geq \kappa_\gamma y_n} \leq \lambda \left( 2e^{-1} ||h||^2_{L^2(d\phi)} ||\rho||^2_{\ell^\tau} \kappa_\gamma^{-2} \right)^{\log(n)} = \OO{n^{-\gamma}},
$$
if $\kappa_\gamma^2 > 2  ||h||^2_{L^2(d\phi)} ||\rho||^2_{\ell^\tau} \exp(\gamma-1)$.

%% alpha tau=1
\noindent\underline{\bf Case  $\bm{ \alpha \tau =1}$ :}
%if we assume that $L(\cdot)$ is an increasing function, it comes that there exists a constant $\lambda>0$
%\begin{equation} \label{cas2eq1}
%\sum_{i \leq n} |\rho(i)|^\tau  \leq \lambda \log(n) L(n)^\tau  \quad \mbox{ and } \quad \sum_{i \leq n} %|\rho(i)|^j =\mathcal{O}(1), \mbox{ for } j>\tau.
%\end{equation}
Using the proof of Theorem $1^\prime$ of \Citet{Breuer83}, we can prove that for all $q\geq 1$
%\begin{eqnarray}
%\Esp \bigg( \big(\overline{Y}_n \big)^{2q} \bigg) &\leq& \kappa n^{-q} L^\prime(n)^{q}
%\left(\sum_{|i|\leq n} \sum_{j\geq \tau} \frac{(c_j)^2}{j!}\rho(i)^j \right)^q  \nonumber \\
%&\leq &\kappa\lambda^q n^{-q} \log(n)^{q} L(n)^{\tau q} \left(\frac{(c_\tau)^2}{\tau!} \right)^q \frac{(2q)!}{2^q q!}. \label{cas2eq2}
%\end{eqnarray}
\begin{eqnarray}
\Esp \left( \big( n^{1/2} L_\tau(n)^{-1/2} \overline{Y}_n \big)^{2q} \right) &\leq& \lambda \frac{2q!}{2^q q!}
\Esp   \left( \big({n}^{1/2} L_\tau(n)^{-1/2} \overline{Y}_n \big)^{2} \right)^q \nonumber \\
&\leq& \lambda \frac{2q!}{2^q q!} \left(
\sum_{j \geq \tau} \frac{(c_j)^2}{j!} \lim_{n\to +\infty} L_\tau(n)^{-1} \sum_{|i|\leq n}|\rho(i)|^j
\right)^q \nonumber \\
&\leq & \lambda \frac{2q!}{2^q q!} ||h||_{L^2(d\phi)}^{2q} .\label{cas2eq1}
\end{eqnarray}
Then from Chebyshev's inequality, we have for all $q\geq 1$
$$ \Prob{| \overline{Y}_n | \geq \kappa_\gamma y_n} \leq \lambda \frac{L_\tau(n)^q}{n^q y_n^{2q}} \; \frac{2q!}{2^q q!} \left( ||h||_{L^2(d\phi)}^2  \kappa_\gamma^{-2} \right)^q.
$$
From (\ref{defeps}) by choosing $q=[\log(n)]$, we obtain
$$
\Prob{| \overline{Y}_n | \geq \kappa_\gamma y_n} \leq \lambda \left( 2 e^{-1}||h||_{L^2(d\phi)}^{2} \; \kappa_\gamma^{-2} \right)^{\log(n)} = \OO{n^{-\gamma}},
$$
if $ \kappa_\gamma^2 >  2  ||h||_{L^2(d\phi)}^{2}   \times \exp(\gamma-1)$.\\
\underline{\bf Case  $\bm{ \alpha \tau <1}$ :}
Denote by $k_{\alpha}$ the lowest integer satisfying $k_\alpha \alpha >1$, that is $k_{\alpha}=[1/\alpha]+1$, and for $j\geq \tau$ denote by $Z_j$ the following random variable
$$
Z_j = \frac1n \sum_{i=1}^n \frac{c_j}{j!} H_j(Y(i)).
$$
Denote by $\kappa_{1,\gamma}$ and $\kappa_{2,\gamma}$ two positive constants such that $\kappa_\gamma=\max (\kappa_{1,\gamma}, \kappa_{2,\gamma})$. From the triangle inequality,
\begin{equation} \label{cas3eq1}
\Prob{| \overline{Y}_n | \geq \kappa_\gamma y_n} \leq \mathbb{P} \Big(| \overline{Y}_n -\sum_{j=\tau}^{k_\alpha-1} Z_j|\geq \kappa_{1,\gamma} y_n \Big) +  \sum_{j=\tau}^{k_\alpha-1} \Prob{ |Z_j| \geq \kappa_{2,\gamma} y_n}
\end{equation}
Since
$$ \overline{Y}_n -\sum_{j=\tau}^{k_\alpha-1} Z_j = \frac1n \sum_{i=1}^n \sum_{j \geq k_\alpha} \frac{c_j}{j!} H_j(Y(i)) = \frac1n \sum_{i=1}^n h^\prime(Y(i)),
$$
where $h^\prime(\cdot)$ is a function with Hermite rank $k_\alpha$. Applying Lemma~\ref{moyps} in the case $\alpha \tau>1$, it follows that, for all $\gamma>0$, there exists a constant $\kappa_{1,\gamma}$ such that, for $n$ large enough
\begin{equation} \label{cas3t1}
\mathbb{P} \Big(| \overline{Y}_n -\sum_{j=\tau}^{k_\alpha-1} Z_j|\geq \kappa_{1,\gamma} y_n \Big)=\OO{n^{-\gamma}}.
\end{equation}
Now, let $\tau \leq j < k_\alpha$ and $q\geq 1$, from Theorem 3 of \Citet{Taqqu77}, we have
\begin{eqnarray}
\Prob{|Z_j|\geq \kappa_{2,\gamma} y_n} &\leq& \frac{1}{\kappa_{2,\gamma}^{2q} y_n^{2q}} \left( \frac{c_j}{j!} \right)^{2q} n^{-2q} \;  \Esp\left( \sum_{i_1,\ldots,i_{2q}} H_j(Y(i_1)) \ldots H_j(Y(i_{2q}))\right) \nonumber \\
&\leq & \lambda \frac{L(n)^{j q}}{n^{\alpha j q} y_n^{2q}} \left( \frac{c_j}{j!} \kappa_{2,\gamma}^{-1} \right)^{2q} \mu_{2q}, \label{cas3eq2}
\end{eqnarray}
where $\mu_{2q}$ is a constant such that $\mu_{2q} \leq \left(\frac{2}{1-\alpha j}\right)^q \Esp \left( H_j(Y)^{2q}\right)$. It is also proved in \Citet{Taqqu77} (p. 228), that $\Esp \left( H_j(Y)^{2q}\right) \sim (2jq)!/(2^{jq} (jq)!),$ as $q\to+\infty$. Thus, from Stirling's formula, we obtain as $q \to +\infty$
$$
\Prob{|Z_j|\geq y_n} \leq \lambda \;  \frac{L(n)^{(j-\tau) q}}{n^{\alpha(j-\tau)q} } \log(n)^{-\tau q} q^{jq} \left( \frac{2}{1-\alpha j} \left(\frac{c_j}{j!}\right)^2 \left(\frac{2j}{e}\right)^j \kappa_{2,\gamma}^{-1}\right)^q .
$$
By choosing $q=[\log(n)]$, we finally obtain, as $n \to +\infty$
\begin{equation} \label{cas3t2}
\sum_{j=\tau}^{k_\alpha-1} \Prob{ |Z_j| \geq \kappa_{2,\gamma} y_n} \leq \lambda \left( \frac{2}{1-\alpha \tau} \left(\frac{c_\tau}{\tau!}\right)^2 \left(\frac{2\tau}{e}\right)^\tau \kappa_{2,\gamma}^{-2}\right)^{\log(n)} = \OO{n^{-\gamma}},
\end{equation}
if $\kappa_{2,\gamma}^2 > \frac{2}{1-\alpha \tau} \left(\frac{c_\tau}{\tau!}\right)^2 (2\tau)^\tau  \exp(\gamma-\tau)$. From (\ref{cas3eq1}), we get the result by combining~(\ref{cas3t1}) and~(\ref{cas3t2}).\end{proof}

\begin{corollary} \label{corIneg} Under conditions of Lemma~\ref{moyps}, for all $\alpha>0$, $j\geq 1$ and $\gamma>0$, there exists $q=q(\gamma)\geq 1$ and $\zeta_\gamma>0$ such that
\begin{equation} \label{momHermite}
\Esp\left( \left\{
\frac1n \sum_{i=1}^n H_j(Y(i))
\right\}^{2q} \right) \leq \zeta_\gamma n^{-\gamma}.
\end{equation}
\end{corollary}

\begin{proof}
(\ref{contMom}), (\ref{cas2eq1}) and~(\ref{cas3eq2}) imply that there exists $\lambda=\lambda(q)>0$ such that for all $q\geq 1$, we have
\begin{eqnarray}
\Esp\left( \left\{
\frac1n \sum_{i=1}^n H_j(Y(i))
\right\}^{2q} \right) \leq \lambda(q) n^{-q} &=&
\lambda(q) \times \left\{ \begin{array}{ll}
n^{-q} & \mbox{ if } \alpha j >1 \\
L_{\tau_p}(n) n^{-q} & \mbox{ if } \alpha j =1 \\
L(n)^{\alpha j q} n^{-\alpha j q}  &\mbox{ if } \alpha j <1 \\
\end{array} \right. \nonumber \\
&=& \OO{n^{-\gamma}} .
\end{eqnarray}
Indeed, it is sufficient to choose $q$ such that, $q>\gamma$  if $\alpha j \geq1$ and
 $q>\gamma/\alpha j$ if $\alpha j <1$. \\
\end{proof}

\begin{lemma} \label{lemme-coeff} %Let $u \in \RR$, denote by $g(\cdot)$ a function satisfying Assumption $\bm{A_3(u)}$. Denote by $c_j(u)$ the $j$-th Hermite coefficient of $h_u(\cdot)$, where $h_u(\cdot)$ is defined by~(\ref{defh}) and by $\tau$ its Hermite rank ($\tau \geq 1$). \\
%$(i)$ For all $j \geq \tau$ we have
%\begin{equation} \label{cju}
%c_j(u) = \left\{ \begin{array}{ll}
%\sum_{i=1}^L -\left( \phi(b_i(u)) - \phi(a_i(u)) \right) & \mbox{if } j=1 ,\\
%\sum_{i=1}^L (-1)^j \left( \phi^{(j-1)}(b_i(u)) - \phi^{(j-1)}(a_i(u)) \right)  & \mbox{if } j>1 ,
%\end{array} \right.
%\end{equation}
%where $g_i(\cdot)$ is the restriction of $g(\cdot)$ on $U_i$ and where $a_i(u)< b_i(u)$ are defined by
%$$
%[a_i(u), b_i(u)] = \left\{ \begin{array}{ll}
%\left[ \inf{g_i^{-1}(U_i)} , g_i^{-1}(u) \right] & \mbox{if } g_i \mbox{ is an increasing function},\\
%\left[ g_i^{-1}(u), \sup{g_i^{-1}(U_i)}  \right] & \mbox{otherwise.} \\
%\end{array} \right.
%$$
Let $0<p<1$, denote by $g(\cdot)$ a function satisfying Assumption $\bm{A_4(\qt{}{p})}$ and by $(x_n)_{n\geq 1}$ a sequence with real components, such that $x_n \to 0$, as $n \to +\infty$. Then, for all $j \geq 1$, there exists a positive constant $d_j=d_j(\qt{}{p})<+\infty$ such that, for $n$ large enough
\begin{equation} \label{cjun}
\left| c_j(\qt{}{p}+x_n) - c_j(\qt{}{p}) \right| \leq d_j \; |x_n| .
\end{equation}
\end{lemma}

\begin{proof} Let $j \geq 1$, under Assumption $\bm{A_4(\qt{}{p})}$, for $n$ large enough, $\qt{}{p}+x_n \in \cup_{i=1}^L g(U_i)$. Thus, for $n$ large enough,
\begin{eqnarray}
c_j(\qt{}{p}+x_n) -c_j(\qt{}{p})&=& \int_{\RR} \left( h_{\qt{}{p}+x_n}(t) -h_{\qt{}{p}}(t)  \right) H_j(t) \phi(t) dt \nonumber \\
&=& \sum_{i=1}^L \int_{U_i} \left( \bm{1}_{g_i(t)\leq \qt{}{p}+x_n} - \bm{1}_{g_i(t)\leq \qt{}{p}}\right) H_j(t)\phi(t) dt \nonumber \\
&=& \sum_{i=1}^L \int_{m_{i,n}}^{M_{i,n}} (-1)^j \phi^{(j)}(t)dt, \nonumber \\
&=& \left\{ \begin{array}{ll}
\sum_{i=1}^L -\left( \phi(M_{i,n}) - \phi(m_{i,n}) \right) & \mbox{if } j=1 ,\\
\sum_{i=1}^L (-1)^j \left( \phi^{(j-1)}(M_{i,n}) - \phi^{(j-1)}(m_{i,n}) \right)  & \mbox{if } j>1 ,
\end{array} \right. \nonumber
\end{eqnarray}
where $g_i(\cdot)$ is the restriction of $g(\cdot)$ to $U_i$, and where $m_{i,n}$ (resp. $M_{i,n}$) is the minimum (resp. maximum) between $g_i^{-1}(\qt{}{p}+x_n)$ and $g_i^{-1}(\qt{}{p})$.
We leave the reader to check that there exists a positive constant $d_j$, such that, for $n$ large enough
$$\left| c_j(\qt{}{p}+x_n) - c_j(\qt{}{p}) \right| \leq  \! d_j  |x_n| \times
\left\{  \begin{array}{ll}
\!\! \sum_{i=1}^L \left| \phi^{(j)} ( g_i^{(-1)}(u) ) \;  \big(g_i^{(-1)} \big)^\prime(u) \right| \!\! & \mbox{if } j=1,2  \\
\!\! \sum_{i=1}^L \left| \phi^{(j-2)}( g_i^{(-1)}(u)) \;  \big(g_i^{(-1)} \big)^\prime(u)\right|  \!\! & \mbox{if } j>2,
\end{array} \right. $$
which is the desired result.
\end{proof}

%\begin{remark} \label{rem-diffCoeff}
%Note that (\ref{cjun}) is valid for all $j\geq 1$ which means that for $1\leq j < \tau_p$, we have
%\begin{equation} \label{cjBis}
%c_j(\qt{}{p}+x_n) \leq d_j | x_n|.
%\end{equation}
%\end{remark}

\begin{lemma} \label{lemme-controleQt}
Under conditions of Theorem~\ref{bahadur}, there exists a constant denoted by $\kappa_\varepsilon= \kappa_\varepsilon(\alpha,\tau_p)$, such that, we have almost surely, as $n\to +\infty$
\begin{equation} \label{controleQt}
\left| \qtEst{}{p}{g(Y)} - \qt{g(Y)}{p}  \right| \leq \varepsilon_n,
\end{equation}
where $\varepsilon_n= \varepsilon_n(\alpha,\tau(\qt{}{p}))= \kappa_{\varepsilon} y_n(\alpha,\tau(\qt{}{p})$, $y_n(\cdot,\cdot)$ being defined by (\ref{defeps}).
\end{lemma}

\begin{proof} We have
\begin{equation} \label{qteq1}
\Prob{ \left| \qtEstSimp{p} - \qt{}{p}  \right| \geq \varepsilon_n } = \Prob{ \qtEstSimp{p} \leq \qt{}{p} - \varepsilon_n} + \Prob{ \qtEstSimp{p} \geq  \qt{}{p}+ \varepsilon_n }.
\end{equation}
Using Lemma 1.1.4 $(iii)$ of \Citet{Serfling80}, we have
\begin{equation} \label{qteq2}
\Prob{ \qtEstSimp{p} \leq \qt{}{p} - \varepsilon_n} \leq  \Prob{ \cdfEstSimp{\qt{}{p}-\varepsilon_n} \geq p}.
\end{equation}
Under Assumption $\bm{A_4(\qt{}{p})}$, for $n$ large enough
$$
p -\cdf{}{\qt{}{p}-\varepsilon_n} = \pdf{}{\qt{}{p}} \varepsilon_n + \oo{\varepsilon_n} \geq \frac{\pdf{}{\qt{}{p}}}{2} \varepsilon_n.
$$
Consequently, for $n$ large enough and from~(\ref{qteq2})
\begin{equation}\label{qteq4}
\Prob{ \qtEstSimp{p} \leq \qt{}{p} - \varepsilon_n} \leq \Prob{ \cdfEstSimp{\qt{}{p}-\varepsilon_n} -\cdf{}{\qt{}{p}-\varepsilon_n} \geq   \frac{\pdf{}{\qt{}{p}}}{2} \varepsilon_n}.
\end{equation}
Define $\tau_{p,n}=\tau(\qt{}{p}-\varepsilon_n)$, from Lemma~\ref{lemme-coeff}, we have for $n$ large enough
\begin{equation} \label{qteq5}
\cdfEstSimp{\qt{}{p}-\varepsilon_n} - \cdf{}{\qt{}{p}-\varepsilon_n} \geq 2 \left(\cdfEstSimp{\qt{}{p}} - \cdf{}{\qt{}{p}}\right)  + 2\varepsilon_n \sum_{j \in J_n} Z_{n,j},
\end{equation}
where
$$
J_n  =\left\{ \begin{array}{ll}
\{\tau_p < j \leq \tau_{p,n}\}  & \mbox{ if } \tau_{p,n}>\tau_p, \\
\emptyset & \mbox{ if } \tau_{p,n}=\tau_p ,\\
\{\tau_{p,n} \leq j < \tau_{p}\}  & \mbox{ if } \tau_{p,n}<\tau_p. \\
\end{array} \right.
 \quad \mbox{ and } \quad Z_{n,j}= \frac1n\sum_{i=1}^n \frac{d_j}{j!} H_j(Y(i)).
$$
Now, define $c_\varepsilon=\kappa_\varepsilon \pdf{}{\qt{}{p}}/4$. Let $\gamma>0$, (\ref{momHermite}) implies that there exists $q\geq 1$ such that, for $n$ large enough
\begin{eqnarray}
\Prob{ | 2\varepsilon_n Z_n| \geq \frac{\pdf{}{\qt{}{p}}}{2} \varepsilon_n} &\leq & \sum_{j\in J_n}
\Prob{|Z_{n,j}|>c_\varepsilon} \nonumber \\
&\leq& \sum_{j\in J_n}  \frac1{c_\varepsilon^{2q}} \Esp\left( Z_{n,j}^{2q}\right)  =\OO{n^{-\gamma}}. \label{termeZn}
\end{eqnarray}
Let us fix $\gamma=2$. From~(\ref{qteq4}), (\ref{qteq5}) and (\ref{termeZn}) and from Lemma~\ref{moyps} (applied to the function $h_{\qt{}{p}}(\cdot)$), we obtain
$$
\Prob{ \qtEstSimp{p} \leq \qt{}{p} - \varepsilon_n} \leq \Prob{  |\cdfEstSimp{\qt{}{p}} - \cdf{}{\qt{}{p}}| \geq c_\varepsilon \varepsilon_n} + \OO{n^{-2}}= \OO{n^{-2}},
$$
if $c_\varepsilon>\kappa_2$ that is if $\kappa_\varepsilon > 4/\pdf{}{\qt{}{p}} \kappa_2$.

Let us now focus on the second right-hand term of~(\ref{qteq1}). Following the sketch of this proof, we may also obtain, for $n$ large enough
$$\Prob{ \qtEstSimp{p} \geq  \qt{}{p}+ \varepsilon_n } = \OO{n^{-2}},
$$
if $\kappa_\varepsilon > 4/\pdf{}{\qt{}{p}} \kappa_2$. Thus, for $n$ large enough $\Prob{ \left|\qtEstSimp{p} - \qt{}{p} \right| \geq \varepsilon_n } = \OO{n^{-2}}$, which leads to the result thanks to Borel-Cantelli's Lemma.
\end{proof}

The following Lemma is an analogous result obtained by Bahadur in the i.i.d. framework, see Lemma~E p.97 of \Citet{Serfling80}.

\begin{lemma} \label{lemme-bahadur}
Under conditions of Theorem~\ref{bahadur}, denote by  $\Delta(z)$ for $z \in\RR$ the random variable,
$\Delta(z) = \cdfEst{}{z}{g(Y)} - \cdf{g(Y)}{ z}$.
Then, we have almost surely, as $n \to +\infty$
\begin{equation} \label{controleSup}
S_n(\qt{g(Y)}{p},\varepsilon_n(\alpha,\tau_p)) = \sup_{|x|\leq \varepsilon_n}
\left|  \Delta(\qt{g(Y)}{p}+x) -  \Delta(\qt{g(Y)}{p})\right|
= \Oas{ r_n(\alpha,\tauB_p)},
\end{equation}
where $\varepsilon_n=\varepsilon_n(\alpha,\tau_p)$ is defined by (\ref{controleQt}) and $r_n(\alpha,\tauB_p)$ is defined by~(\ref{defrnAlpha}).
\end{lemma}

\begin{proof}
Put $\varepsilon_n=\varepsilon_n(\alpha,\tau_p)$ and $r_n=r_n(\alpha,\tauB_p)$.
Denote by $(\beta_n)_{n\geq 1}$ and $(\eta_{b,n})_{n\geq 1}$ the following two sequences
$$\beta_n=\big[n^{3/4}\varepsilon_n \big]\quad \mbox{ and } \quad \eta_{b,n}=\qt{}{p} + \varepsilon_n \frac{b}{\beta_n},$$
for $b=-\beta_n,\ldots,\beta_n$. Using the monotonicity of $\cdf{}{\cdot}$ and $\cdfEstSimp{\cdot}$, we have,
\begin{equation} \label{compSn}
S_n(\qt{}{p},\varepsilon_n) \leq \max_{-\beta_n\leq b \leq \beta_n} |M_{b,n}| + G_n,
\end{equation}
where  $M_{b,n}=\Delta(\eta_{b,n})- \Delta(\qt{}{p})$ and $G_n = \max_{-\beta_n\leq b \leq \beta_n-1} \left( \cdf{}{\eta_{b+1,n}} -  \cdf{}{\eta_{b,n}} \right).$
Under Assumption $\bm{A_4(\qt{}{p})}$, we have for $n$ large enough
\begin{equation} \label{compGn}
G_n \leq \left( \eta_{b+1,n}-\eta_{b,n} \right) \times \sup_{|x|\leq \varepsilon_n} \pdf{}{\qt{}{p}+x}  = \OO{n^{-3/4}}.
\end{equation}
The proof is finished if one can prove that for all $\gamma>0$ (in particular $\gamma=2$) and for all $b$,  there exists $\kappa^\prime_\gamma$ such that
\begin{equation} \label{objectif}
\Prob{ |M_{b,n}| \geq \kappa_\gamma^\prime r_n} = \OO{n^{-\gamma}}.
\end{equation}
Indeed, since $\beta_n=\OO{n^{1/2+\delta}}$ for all $\delta>0$, if (\ref{objectif}) is true, then we have
\begin{eqnarray}
\mathbb{P} \big( \max_{-\beta_n\leq b \leq \beta_n} |M_{b,n}| \geq \kappa_2^\prime r_n(\alpha,\tau_p)\big) \; 
&\leq& \; (2\beta_n+1) \times \max_{-\beta_n\leq b \leq \beta_n}\Prob{ |M_{b,n}| \geq \kappa_2^\prime r_n} \nonumber\\  
&=& \OO{n^{-3/2+\delta}}. \nonumber
\end{eqnarray}
Thus, from Borel-Cantelli's Lemma, we have, almost surely
$$\max_{-\beta_n\leq b \leq \beta_n} |M_{b,n}| =\Oas{r_n}$$
And so, from~(\ref{compSn}) and~(\ref{compGn}).
\begin{equation} \label{compSnFinal}
S_n(\qt{}{p},\varepsilon_n)= \Oas{r_n} + \OO{n^{-3/4}}  = \Oas{r_n} ,
\end{equation}
which is the stated result.

So, the rest of the proof is devoted to prove~(\ref{objectif}). For the sake of simplicity, denote by $h_n^\prime(\cdot)$ the function $h_{\eta_{b,n}}(\cdot) - h_{\qt{}{p}}(\cdot)$. For $n$ large enough, the Hermite rank of $h_n^\prime(\cdot)$ is at least equal to $\tauB_p$, that is defined by~(\ref{tauBarre}). In the sequel, we need the following bound for $||h_n^\prime||^2_{L^2(d\phi)}$
$$
||h_n^\prime||^2_{L^2(d\phi)}= \Esp(h_n^\prime(Y)^2) = \omega_n ( 1-\omega_n) \quad \mbox{ with } \quad \omega_n= \big| \cdf{g(Y)}{\eta_{b,n}}- \cdf{g(Y)}{\qt{}{p}} \big|.
$$
As previously, we have $\omega_n= \mathcal{O}(\varepsilon_n)$ and so, there exists $\zeta>0$, such that
\begin{equation} \label{compEh2}
 ||h_n^\prime||^2_{L^2(d\phi)} \leq \zeta \varepsilon_n.
\end{equation}
From now on, in order to simplify the proof, we use the following upper-bound
$$
\varepsilon_n=\varepsilon_n(\alpha,\tau_p) \leq \varepsilon_n(\alpha ,\tauB_p),
$$
and with a slight abuse, we still denote $\varepsilon_n=\varepsilon_n(\alpha ,\tauB_p)$.
Note also, that from Lemma~\ref{lemme-coeff}, the $j$-th Hermite coefficient, for some $j \geq \tauB_p$, is given by $c_j(\eta_{b,n})-c_j(\qt{}{p})$. And there exists a positive constant $d_j=d_j(\qt{}{p})$ such that for $n$ large enough
\begin{equation} \label{diffcj}
| c_j(\eta_{b,n})-c_j(\qt{}{p} | \leq d_j \; \varepsilon_n  \frac{|b|}{\beta_n} \; \leq \;  d_j \; \varepsilon_n .
\end{equation}
%Note also, that for all $1\leq j \leq \tauB_p$, we have $|c_j(\eta_{b,n})|\leq d_j \; \varepsilon_n$.
We now proceed like in the proof of Lemma~\ref{moyps}.

\noindent \underline{\bf Case $\bm{\alpha \tauB_p>1}$}: using Theorem~1 of \Citet{Breuer83} and~(\ref{cas1eq1}), we can obtain for all $q\geq 1$
\begin{equation} \label{Mncas1eq1}
\Prob{ |M_{b,n}| \geq \kappa_\gamma^\prime r_n} \leq \lambda
\frac{1}{n^q   r_n^{2q}}  \frac{(2q)!}{2^q q!} \frac1{(\kappa_\gamma^\prime)^{2q}} ||h_n^\prime||^{2q}_{L^2(d\phi)} ||\rho||_{\ell^{\tauB_p}}^{2q} .
\end{equation}
As $q\to +\infty$, we get
$$
\Prob{ |M_{b,n}| \geq \kappa_\gamma^\prime r_n} \leq \lambda \frac{\varepsilon_n^{q}}{n^q   r_n^{2q}} q^q \left( 2\zeta  e^{-1} ||\rho||_{\ell^{\tauB_p}}^{2} \frac1{(\kappa_\gamma^\prime)^{2}}\right)^q.
$$
From~(\ref{defrnAlpha}), (\ref{defeps}) (with $\tau=\tauB_p$) and by choosing $q=[\log(n)]$, we have
\begin{equation} \label{Mncas1eq2}
\Prob{ |M_{b,n}| \geq \kappa^\prime_\gamma  r_n} \leq \lambda \left( 2\zeta \kappa_\varepsilon e^{-1} ||\rho||_{\ell^{\tauB_p}}^{2} \frac1{(\kappa_\gamma^\prime)^{2}}\right)^{\log(n)}  = \OO{n^{-\gamma}},
\end{equation}
if ${\kappa_\gamma^\prime}^{2} > 2\zeta \kappa_\varepsilon  ||\rho||_{\ell^{\tau_p}}^{2} \exp(\gamma-1)$.

\noindent \underline{\bf Case $\bm{\alpha \tauB_p=1}$}
from~(\ref{cas2eq1}), we can obtain  for all $q\geq 1$
\begin{eqnarray*}
\Esp \left( M_{b,n}^{2q} \right) &\leq& \lambda \frac{(2q)!}{2^q q!} \; \frac{L_{\tauB_p}(n)^q}{n^q}  ||h_n^\prime||^{2q}_{L^2(d\phi)} \; \leq \; \lambda \; \zeta^q \;\frac{(2q)!}{2^q q!} \; \frac{L_{\tauB_p}(n)^q \varepsilon_n^q}{n^q} \nonumber \\
&\leq & \lambda \; \frac{L_{\tauB_p}(n)^q \varepsilon_n^q}{n^q} \; (2\zeta e^{-1})^q  \; q^q.
\end{eqnarray*}
From~(\ref{defrnAlpha}), (\ref{defeps}) (with $\tau=\tauB_p$), by choosing $q=[\log(n)]$, we have
\begin{eqnarray}
\Prob{ |M_{b,n}| \geq \kappa_\gamma^\prime r_n}  &\leq&  \frac{1}{{\kappa_\gamma^\prime}^{2q} r_n^{2q}}
\Esp \left( M_{b,n}^{2q} \right) \nonumber \\
&\leq& \lambda \left(
2 \zeta \; \kappa_{\varepsilon} \; e^{-1}\frac{d_{\tau_p}^2}{\tau_p!} \frac1{ {\kappa_\gamma^\prime}^2 }
\right)^{\log(n)}  = \OO{n^{-\gamma}}, \nonumber
\end{eqnarray}
if ${\kappa_\gamma^\prime}^2 > 2 \zeta \kappa_\varepsilon d_{\tau_p}^2 /\tau_p! \exp(\gamma-1)$.

\noindent \underline{\bf Case $\bm{\alpha \tauB_p<1}$}:
denote by $(r_{1,n})_{n\geq 1}$ and by $(r_{2,n})_{n\geq 1}$ the following two sequences
\begin{equation} \label{def2rn}
r_{1,n}= n^{-1/2-\alpha\tauB_p/4} \log(n)^{\tauB_p/4+1/2} L(n)^{\tauB_p/4} \; \mbox{ and } \;
r_{2,n}= n^{-\alpha \tauB_p} \log(n)^{\tauB_p} L(n)^{\tauB_p}.
\end{equation}
Note that $\max \left( r_{1,n}, r_{2,n} \right)$ is equal to  $r_{1,n}$, when $2/3< \alpha\tauB_p<1$ and to $r_{2,n}$, when $0< \alpha\tauB_p \leq 2/3$. So, in order to obtain~(\ref{objectif}) in the case $0<\alpha\tauB_p<1$, it is sufficient to prove that there exists $\kappa_\gamma^\prime$ such that, for $n$ large enough
$$
\Prob{ |M_{b,n}| \geq \kappa_\gamma^\prime \max ( r_{1,n}, r_{2,n} )} = \OO{n^{-\gamma}}.
$$
Denote by $k_\alpha$ the integer $[1/\alpha]+1$ for which $\alpha k_\alpha>1$, and by $Z_{j,n}$ for $\tauB_p\leq j < k_\alpha$ the random variable defined by
$$
Z_{j,n} = \frac1n \sum_{i=1}^n \frac{c_j(\eta_{b,n})-c_j(\qt{}{p})}{j!} \; H_j(Y(i)).
$$
From the triangle inequality, we have
\begin{equation}\label{Mncas3eq1}
\Prob{ |M_{b,n}| \geq \kappa_\gamma^\prime \max ( r_{1,n}, r_{2,n} )} \leq \mathbb{P} \big( \big| M_{b,n} -\sum_{j=\tauB_p}^{k_\alpha-1}Z_{j,n} \big| \geq \kappa_\gamma^\prime r_{1,n} \big) + \sum_{j=\tauB_p}^{k_\alpha-1} \Prob{|Z_{j,n}|\geq \kappa_\gamma^\prime r_{2,n}} .
\end{equation}
Since,
$$  M_{b,n} -\sum_{j=\tau_p}^{k_\alpha-1} Z_{j,n}= \frac1n \sum_{i=1}^n \sum_{j \geq k_\alpha} \frac{c_j(\eta_{b,n})-c_j(\qt{}{p})}{j!} H_j(Y(i)) = \frac1n\sum_{i=1}^n  h_n^{\prime\prime}(Y(i)),
$$
where $h_n^{\prime\prime}(\cdot)$ is a function with Hermite rank $k_\alpha$, such that $\alpha k_\alpha>1$, we have from~(\ref{Mncas1eq1})
\begin{equation}\label{Mncas3eq2}
\mathbb{P} \big( \big| M_{b,n} -\sum_{j=\tauB_p}^{k_\alpha-1} Z_{j,n}\big| \geq \kappa_\gamma^\prime r_{1,n} \big) \leq  \lambda \; \frac1{n^q r_{1,n}^{2q}} ||h_n^\prime||_{L^2(d\phi)}^{2q} \frac{(2q)!}{2^q q!}
\frac{1}{{\kappa_\gamma^\prime}^{2q}} ||\rho||_{\ell^{k_\alpha}}^{2q}
\end{equation}
for all $q\geq 1$. From~(\ref{compEh2}), we obtain, as $q\to +\infty$
$$
\mathbb{P} \big( \big| M_{b,n} -\sum_{j=\tau_p}^{k_\alpha-1} Z_{j,n}\big| \geq \kappa_\gamma^\prime r_{1,n} \big) \leq  \lambda \frac{\varepsilon_n^q}{n^q r_{1,n}^{2q}}  q^q   \left( 2\zeta e^{-1} ||\rho||_{\ell^{k_\alpha}}^{2} {\kappa_\gamma^{\prime}}^{-2} \right)^q .
$$
From~(\ref{defeps}) (with $\tau=\tauB_p$), (\ref{def2rn}) and by choosing $q=[\log(n)]$, we obtain
\begin{equation} \label{t1cas3eq1}
\mathbb{P} \big( \big| M_{b,n} -\sum_{j=\tau_p}^{k_\alpha-1} Z_{j,n}\big| \geq \kappa_\gamma^\prime r_{1,n} \big) \leq \lambda  \left( 2\zeta e^{-1} ||\rho||_{\ell^{k_\alpha}}^{2} \kappa_\varepsilon {\kappa_\gamma^{\prime}}^{-2} \right)^{\log(n)}  = \OO{n^{-\gamma}},
\end{equation}
if ${\kappa_\gamma^\prime}^2 > \kappa_{1,\gamma}^\prime= 2\zeta ||\rho||_{\ell^{k_\alpha}}^{2} \kappa_\varepsilon \exp(\gamma-1)$. Now, concerning the last term of~(\ref{Mncas3eq1}), from~(\ref{cas3eq2}), we can prove, for all $\tauB_p\leq j < k_\alpha$
$$
\Prob{Z_{j,n} \geq \kappa_\gamma^\prime r_{2,n} } \leq
\lambda \frac{L(n)^{j q}}{n^{\alpha j q} \; r_{2,n}^{2q}}
\; \frac1{ {\kappa_\gamma^\prime}^{2q} }\;
\left( \frac{
c_j(\eta_{b,n})-c_j(\qt{}{p}) }{j!} \right)^{2q} \mu_{2q},
$$
where $\mu_{2q}$ is a constant such that, as $q\to +\infty$,
$$\mu_{2q} \leq \lambda \left(\frac{2}{1-\alpha j}\right)^q  \frac{(2jq)!}{2^{jq} (jq)!}.$$
From~(\ref{diffcj}), we have, as $q\to+\infty$
$$
\Prob{Z_{j,n} \geq \kappa_\gamma^\prime r_{2,n} } \leq
\lambda  \frac{\varepsilon_n^{2q} L(n)^{j q}}{n^{\alpha j q} \; r_{2,n}^{2q}} q^{jq} \left(
\frac{2}{1-\alpha j} \left(\frac{2j}{e} \right)^j \; d_j^2 \; {\kappa_\gamma^\prime}^{-2}
\right)^{2q}.
$$
From~(\ref{defrnAlpha}), (\ref{defeps}) (with $\tau=\tauB_p$) by choosing $q=[\log(n)]$, we have, as $n \to +\infty$
$$
\Prob{Z_{j,n} \geq \kappa_\gamma^\prime r_{2,n} } \leq
\lambda \left(\frac{  \log(n) L(n)}{ n^{\alpha}} \right)^{(j-\tauB_p)q} \; \left( \frac{2}{1-\alpha j} \left(\frac{2j}{e} \right)^j \; d_j^2 \;{\kappa_\varepsilon^2}  \; {\kappa_\gamma^\prime}^{-2} \right)^{q}.
$$
Consequently, as $n \to +\infty$, we finally obtain
\begin{equation} \label{t2cas3eq1}
\sum_{j=\tauB_p}^{k_\alpha-1} \Prob{Z_{j,n} \geq \kappa_\gamma^\prime r_{2,n} } \leq
\lambda \left( \frac{2}{1-\alpha \tauB} \left(\frac{2\tauB}{e} \right)^{\tauB} \; d_{\tauB}^2 \;\kappa_\varepsilon^2  \; {\kappa_\gamma^\prime}^{-2} \right)^{\log(n)} = \OO{n^{-\gamma}},
\end{equation}
if ${\kappa_\gamma^\prime}^{2} >  \kappa_{2,\gamma}^\prime= \frac{2}{1-\alpha \tauB} \left(\frac{2\tauB}{e} \right)^{\tauB} \; d_{\tauB}^2 \;\kappa_\varepsilon^2 \exp(\gamma-\tauB)$. Let us choose $\kappa_\gamma^\prime$ such that ${\kappa_\gamma^\prime}^2 > \max( \kappa_{1,\gamma}^\prime, \kappa_{2,\gamma}^\prime)$. Then, by combining~(\ref{t1cas3eq1}) and~(\ref{t2cas3eq1}), we deduce from~(\ref{Mncas3eq1}) that, for every $\gamma>0$
$$
\Prob{ |M_{b,n}| \geq \kappa_\gamma^\prime \max ( r_{1,n}, r_{2,n} )} = \OO{n^{-\gamma}} ,
$$
and so, (\ref{objectif}) is proved.
\end{proof}

\subsection{Proof of Theorem~\ref{bahadur}}

\begin{proof}
Let us detail the proof presented in Section~\ref{intro-proof}. We have
$$
\frac{p- \cdfEstSimp{\qt{}{p}}}{f(\qt{}{p})} -\left( \qtEstSimp{p}-\qt{}{p}\right)  = A(p) + B(p) + C(p)$$
with $A(p)$, $B(p)$ and $C(p)$ respectively defined by~(\ref{defA}),~(\ref{defB}) and~(\ref{defC}).
Under Assumption $\bm{A_4(\qt{}{p})}$, from Lemma~\ref{lemme-controleQt} and Taylor's theorem  we have almost surely, as $n\to +\infty$
$$
C (p) \leq \sup_{|x|\leq \varepsilon_n(\alpha,\tau_p)} F_{g(Y)}^{\prime\prime}(\qt{}{p}+x) \left( \qtEstSimp{p}-\qt{}{p}\right)^2 =
\Oas{\varepsilon_n(\alpha,\tau_p)^2}.
$$
From the definition of sample quantile, we have almost surely, see {\it e.g.} \Citet{Serfling80}, $A(p)= \Oas{n^{-1}}$.
Now, by combining  Lemma~\ref{lemme-controleQt} and Lemma~\ref{lemme-bahadur}, we have almost surely $B(p)=\Oas{r_n(\alpha,\tauB_p)}$.
Thus, we finally obtain
$$
\qtEstSimp{p} - \qt{}{p} = \frac{ p - \cdfEstSimp{\qt{}{p}}  }{ \pdf{}{\qt{}{p}} } +  \Oas{n^{-1}}+\Oas{r_n(\alpha,\tauB_p)} + \Oas{\varepsilon_n(\alpha,\tau_p)^2},
$$
which leads to the result by noticing that $\varepsilon_n(\alpha,\tau_p)^2=\OO{r_n(\alpha,\tauB_p)}$.
\end{proof}

%\newpage

\subsection{Auxiliary Lemmas for the proof of Theorem~\ref{convHest}}

Let $0<p_0 \leq p_1<1$.

\begin{lemma} \label{lemme-supqQuant}
Under conditions of Theorem~\ref{unifBahadur}, there exists a constant denoted by $\theta=\theta(\alpha,\tau_{p_0,p_1})$ such that, we have almost surely, as $n\to +\infty$
\begin{equation} \label{eq-supQuant}
T= \sup_{p_0\leq p \leq p_1} \left| \qtEst{}{p}{g(Y)}-\qt{g(Y)}{p}\right| \leq \varepsilon_n(\alpha,\tau_{p_0,p_1}),
\end{equation}
where $\varepsilon_n = \varepsilon_n(\alpha,\tau_{p_0,p_1})= \theta y_n(\alpha,\tau_{p_0,p_1})$ and $y_n$ is given by (50).
\end{lemma}

\begin{proof}
Define $p_{j,n}=p_0+ \frac j{[n^{3/2}]} (p_1-p_0)$ for $j=0,\ldots,[n^{3/2}]$, and let $p \in [p_0,p_1]$. Using the monotonicity of $\widehat{\xi}(\cdot)$ and $\xi(\cdot)$, there exists some $j$ such that $p\in [p_{j,n},p_{j+1,n}]$ and such that
\begin{eqnarray*}
\qtEstSimp{p}-\qt{}{p}&\leq & \qtEstSimp{p}-\qtEstSimp{p_{j+1,n})} +\qtEstSimp{p_{j+1,n})}-\qt{}{p} \\
&\leq& \qtEstSimp{p_{j+1,n}}-\qt{}{p_{j+1,n})} + \qt{}{p_{j+1,n})}-\qt{}{p_{j,n})} + \qt{}{p_{j,n})}-\qt{}{p}\\
&\leq& \qtEstSimp{p_{j+1,n}}-\qt{}{p_{j+1,n})} + \qt{}{p_{j+1,n})}-\qt{}{p_{j,n})}.
\end{eqnarray*}
This leads to
\begin{equation} \label{eqT}
T \leq \max_{j=0,\ldots,[n^{3/2}]}\left|\qtEstSimp{p_{j,n}}-\qt{}{p_{j,n})} \right| +
\max_{j=0,\ldots,[n^{3/2}]-1}  \left|\qt{}{p_{j+1,n}}-\qt{}{p_{j,n})} \right|.
\end{equation}
Under Assumption $\bm{A_5(p_0,p_1)}$, it comes
\begin{equation} \label{max32}
\max_{j=0,\ldots,[n^{3/2}]-1}  \left|\qt{}{p_{j+1,n}}-\qt{}{p_{j,n})} \right|= \OO{n^{-3/2}}.
\end{equation}
Now, following the proof of Lemma~\ref{lemme-controleQt}, one can prove that there exists some constant $\theta(\alpha,\tau_{p_0,p_1})$ such that for all $j=0,\ldots,[n^{3/2}]$,
$$\Prob{|\qtEstSimp{p_{j,n}}-\qt{}{p_{j,n}}|\geq \theta y_n(\alpha,\tau_{p_0,p_1})}=\OO{n^{-3}}.$$
Therefore, as $n \to +\infty$,
\begin{eqnarray*}
\Prob{\max_{j=0,\ldots,[n^{3/2}]}|\qtEstSimp{p_{j,n}}-\qt{}{p_{j,n}}|\geq \varepsilon_n}
&\leq& ([n^{3/2}]+1) \max_{j=0,\ldots,[n^{3/2}]} \Prob{|\qtEstSimp{p_{j,n}}-\qt{}{p_{j,n}}|\geq \varepsilon_n} \\
&=& \OO{n^{-3/2}}.
\end{eqnarray*}
which, combined with (\ref{eqT}), (\ref{max32}) and Borel-Cantelli's Lemma, leads to the result.
\end{proof}

The following result is an extension of Lemma~\ref{lemme-bahadur} and Theorem~4.2 obtained by \Citet{Sen71}.

\begin{lemma} \label{lemme-unifBahadur}
Under Assumptions of Theorem~\ref{unifBahadur} and following Lemma~\ref{lemme-bahadur}, we have almost surely, as $n \to +\infty$
\begin{equation} \label{eq-lemmeUnif}
S_n^\star = \sup_{
\begin{array}{c}
x,y\in [\qt{}{p_0},\qt{}{p_1}] \\
|x-y|\leq \varepsilon_n(\alpha,\tau_{p_0,p_1})
\end{array}
} \left| \Delta(x)-\Delta(y)\right| = \Oas{ r_n(\alpha,\tau_{p_0,p_1})}
\end{equation}
where $\tau_{p_0,p_1}$ is defined by (\ref{taup0p1}).
\end{lemma}

\begin{proof}
Set $\varepsilon_n=\varepsilon_n(\alpha,\tau_{p_0,p_1})$ and  $r_n=r_n(\alpha,\tau_{p_0,p_1})$. Define $\xi_{j,n}= \qt{}{p_0} + \frac{j}{p_n}(\qt{}{p_1}-\qt{}{p_0})$   for $j=0,\ldots,p_n$ with $p_n=\left[\varepsilon_n^{-1}\right]$, and let $x,y \in [\qt{}{p_0},\qt{}{p_1}]$ such that $|x-y|\leq \varepsilon_n$. Two cases may occur
\begin{itemize}
\item If there exists some $j$ such that $x,y \in [\xi_{j,n} , \xi_{j+1,n}]$ then
$$
|\Delta(x)-\Delta(y)| \leq |\Delta(x) - \Delta(\xi_{j,n})| + |\Delta(\xi_{j,n})-\Delta(y)| \leq 2 \times S_n(\xi_{j,n},\varepsilon_n)
$$
\item Otherwise and witout loss of generality, there exists $j,k$ with $k>j$ such that $x\in [\xi_{j,n},\xi_{j+1,n}]$ and $y\in [\xi_{k,n},\xi_{k+1,n}]$. Since $|x-y|\leq \varepsilon_n$, it follows that $|\xi_{k,n}-\xi_{j+1,n}|\leq \varepsilon_n$. Then,
\begin{eqnarray*}
|\Delta(x)-\Delta(y)| &\leq& |\Delta(x) - \Delta(\xi_{k,n})|+ |\Delta(\xi_{k,n}) - \Delta(\xi_{j+1,n})| + |\Delta(\xi_{j+1,n}) - \Delta(y)| \\
&\leq & S_n(\xi_{k,n},\varepsilon_n) + 2 \times S_n(\xi_{j+1,n},\varepsilon_n).
\end{eqnarray*}
\end{itemize}
In other words, for all $x,y$ one may obtain
$$
|\Delta(x) - \Delta(y) | \leq 3 \times \max_{0\leq j \leq p_n} S_n(\xi_{j,n},\varepsilon_n).
$$
Hence, $S_n^\star \leq 3 \times \max_{0\leq j \leq p_n} S_n(\xi_{j,n},\varepsilon_n)$. Now, following the proof of Lemma~\ref{lemme-bahadur}, one may prove that there exists some positive constant $\theta_\gamma$ such that for $n$ large enough and for all $j=0,\ldots,p_n$,
$$
\Prob{S_n(\xi_{j,n},\varepsilon_n)\geq \theta_\gamma r_n} = \OO{n^{-\gamma}}.
$$
And in particular for $\gamma=2$, it comes
\begin{eqnarray*}
\Prob{ \max_{0\leq j \leq p_n} S_n(\xi_{j,n},\varepsilon_n) \geq \theta_2 r_n} &\leq& (p_n+1) \max_{j=0,\ldots,p_n} \Prob{S_n(\xi_{j,n},\varepsilon_n)\geq \theta_\gamma r_n} \\
&=& \OO{\frac{p_n}{n^2}} = \OO{n^{-3/2}},
\end{eqnarray*}
whatever the value of $\alpha\tau_{p_0,p_1}$. This leads to the result by using Borel-Cantelli's Lemma.
 \end{proof}
\subsection{Proof of Theorem~\ref{unifBahadur}}

\begin{proof}
We follow the proof of Theorem~\ref{bahadur}. Let $p\in [p_0,p_1]$ and let $\varepsilon_n=\varepsilon_n(\alpha,\tau_{p_0,p_1})$, then
$$
\frac{p- \cdfEstSimp{\qt{}{p}}}{f(\qt{}{p})} -\left( \qtEstSimp{p}-\qt{}{p}\right)  = A(p) + B(p) + C(p)$$
where $A(p),B(p)$ and $C(p)$  are respectively defined by~(\ref{defA}),~(\ref{defB}) and~(\ref{defC}). Similarly to the proof of Theorem~\ref{bahadur}, one may prove that $\sup_{p_0\leq p \leq p_1}A(p)=\Oas{n^{-1}}$. Under Assumption $\bm{A_5(p_0,p_1)}$,
$C(p)\leq \left(\sup_{|x|\leq \varepsilon_n(\alpha,\tau_p)} F^{\prime \prime}(x+\qt{}{p})\right) \frac{ \left(\qtEstSimp{p}-\qt{}{p}\right)^2  }{f(\qt{}{p})}$. Therefore, for $n$ large enough, $C(p) \leq \lambda \left( \sup_{p_0\leq p \leq p_1}  \left(\qtEstSimp{p}-\qt{}{p}\right)\right)^2$. And from Lemma~\ref{lemme-supqQuant}, this leads to
$$
\sup_{p_0\leq p \leq p_1} C(p) = \Oas{\varepsilon_n(\alpha,\tau_{p_0,p_1})^2}.
$$
In addition, using Lemma~\ref{lemme-unifBahadur}, one also has $\sup_{p_0\leq p \leq p_1} B(p)= \Oas{r_n(\alpha,\tau_{p_0,p_1})}$, which ends the proof.
\end{proof}

\subsection{Auxiliary Lemma for the proof of Theorem~\ref{convHest}}

\begin{lemma} \label{hp}
Consider for $0<p<1$ the function $h_p(\cdot)$, given by
\begin{equation} \label{defhp}
h_p(t) = \bm{1}_{\{  |t|\leq \qt{|Y|}{p}  \}}(t) -p,
\end{equation}
that is the function $h_{\qt{g(Y)}{p}}(\cdot)$ with $g(\cdot)=|\cdot|$. Then by denoting  $c_j^{h_p}$ the $j$-th Hermite coefficient of $h_p(\cdot)$, we have for all $j\geq 1$
\begin{equation} \label{coeffhp}
c_0^{h_p}=c_{2j+1}^{h_p}=0 \quad \mbox{ and } \quad c_{2j}^{h_p}= -2 H_{2j-1}(q)\phi(q),
\end{equation}
where $q=\qt{|Y|}{p} =\Phi^{-1}\left( \frac{1+p}{2} \right)$.
\end{lemma}

\begin{proof}
Since $\Prob{|Y|\leq q}=p$ and $h_p(\cdot)$ is even, we have $c_0^{h_p}=c_{2j+1}^{h_p}=0$, for all $j\geq 1$. Now, (\ref{defHermite}) implies
\begin{eqnarray}
c_{2j}^{h_p} &=& \int_{\RR} h_p(t)H_{2j}(t)\phi(t)dt = 2 \times \int_0^{q } H_{2j}(t)\phi(t)dt \nonumber \\
&=& 2 \times \left[ \phi^{(2j-1)}(t) \right]_{0}^{   q   } = 2 \times \left[ - H_{2j-1}(t) \phi(t) \right]_0^{q} \nonumber \\
&=& - 2 H_{2j-1} (q)  \phi(q). \nonumber
\end{eqnarray}
\end{proof}

\begin{remark} \label{remhp}
Let $g(\cdot)=\widetilde{g}(|\cdot|)$, where $\widetilde{g}(\cdot)$ is a strictly increasing function on $\RR^+$, then for all $0<p<1$, we have
$$
\qt{|Y|}{p} = \widetilde{g}^{-1} \left( \qt{g(Y)}{p} \right).
$$
Consequently, the functions $h_{\qt{g(Y)}{p}}(\cdot)$ for $g(\cdot)=|\cdot|$, $g(\cdot)=|\cdot|^\alpha$ and $g(\cdot)=\log|\cdot|$ are strictly identical. And so, their Hermite decomposition is given by~(\ref{coeffhp}) and their Hermite rank is equal to~2.
\end{remark}

\subsection{Proof of Theorem~\ref{convHest}}

\begin{proof}
$(i)$ Define
\begin{equation}
b_{n}=\frac12 \sum_{m=1}^M B_m \log\left( 1+ \delta_n^{\Vect{a^m}}(0)\right),
\end{equation}
where $\delta_n^{\Vect{a^m}}(0)$ is given by (\ref{deltana}). From~(\ref{defepsmalpha}), (\ref{defepsmlog}), and~(\ref{HestMoinsH}), we have almost surely
\begin{eqnarray}
\hspace*{-0.2cm}\widehat{H}^{\alpha}\!\!-\!\!H &=& \sum_{m=1}^M \frac{B_m}{\alpha} \varepsilon_m^\alpha  \nonumber \\
&=&\sum_{m=1}^M \frac{B_m}{\alpha} \log \left( \frac{
\lqtEst{p}{c}{ |Y^{\Vect{a^m}}|^\alpha } }{
\lqt{|Y|^\alpha}{p}{c} } \right)  +\alpha \times b_n \nonumber \\
&=& \sum_{m=1}^M \frac{B_m}{\alpha \lqt{|Y|^\alpha}{p}{c}} \! \left( 
\lqtEst{p}{c}{ |Y^{\Vect{a^m}}|^\alpha }  - \lqt{|Y|^\alpha}{p}{c}  \right)\!(1+\oas{1}) + \alpha \;b_n. \label{eq1-Halpha} 
\end{eqnarray}
and
\begin{eqnarray}
\widehat{H}^{\log}-H &=& \sum_{m=1}^M  B_m \varepsilon_m^{\log} \nonumber \\
&=& \sum_{m=1}^M B_m \left( 
\lqtEst{p}{c}{ \log|Y^{\Vect{a^m}}|}  - \lqt{\log|Y|}{p}{c} \right) 
+ b_n \label{eq1-Hlog}
\end{eqnarray}
Under Assumption $\bm{A_6(\eta)}$, we have 
\begin{equation} \label{contbn}
b_n= \OO{n^{-\eta}}.
\end{equation}
Moreover, let $i,j\geq1$, under Assumption $\bm{A_1(2\nu)}$, we have, from Lemma~\ref{lemKent}
\begin{equation} \label{rappelCorrelation}
\Esp( Y^{\Vect{a^m}}(i) Y^{\Vect{a^m}}(i+j) ) = \rho^{\Vect{a^m}}(j) = \OO{ |j|^{2H-2\nu}}.
\end{equation}
Then, for all $m=1,\ldots,M$ and for all $k=1,\ldots,K$, from Lemma~\ref{lemme-controleQt} and Remark~\ref{remhp}, we obtain, that almost surely
\begin{eqnarray}
\qtEst{}{p_k}{|Y^{a^m}|^\alpha}- \qt{|Y|^\alpha}{p_k} &=& \Oas{y_n(2\nu-2H,\tau_{p_k})}, \nonumber \\
\qtEst{}{p_k}{\log|Y^{a^m}|} - \qt{\log|Y|}{p_k} &=&  \Oas{y_n(2\nu-2H,
\tau_{p_k})}, \nonumber
\end{eqnarray}
where the sequence $y_n(\cdot,\cdot)$ is defined by~(\ref{defeps}) with $L(\cdot)=1$. The result~(\ref{defyn}) is obtained by combining~(\ref{eq1-Halpha}), (\ref{eq1-Hlog}) and~(\ref{contbn}).

$(ii)$ Let us apply Theorem~\ref{bahadur} to the sequence $\Vect{g(Y^{a^m})}$, for some $m=1,\ldots,M$,  with $g(\cdot)=|\cdot|$, $g(\cdot)=|\cdot|^\alpha$ and $g(\cdot)=\log|\cdot|$. For all $k=1,\ldots,K$, we have almost surely
\begin{eqnarray}
\qtEst{}{p_k}{|Y^{a^m}|}  - \qt{|Y|}{p_k} &=&
\frac{p_k - \cdfEst{}{\qt{|Y|}{p_k} }{|Y^{a^m}|}}{\pdf{|Y|^\alpha}{\qt{|Y|}{p_k}}}
+ \Oas{r_n}
\nonumber \\
\qtEst{}{p_k}{|Y^{a^m}|^\alpha}  - \qt{|Y|^\alpha}{p_k} &=&
\frac{p_k - \cdfEst{}{\qt{|Y|^\alpha}{p_k} }{|Y^{a^m}|^\alpha}}{\pdf{|Y|^\alpha}{\qt{|Y|^\alpha}{p_k}}}
+ \Oas{r_n}
\nonumber\\
\qtEst{}{p_k}{\log|Y^{a^m}|}  - \qt{\log|Y|}{p_k} &=&
\frac{p_k - \cdfEst{}{\qt{\log|Y|}{p_k} }{\log|Y^{a^m}|}}{\pdf{\log|Y|}{\qt{\log|Y|}{p_k}}}
+ \Oas{r_n}, \nonumber
\end{eqnarray}
where, for the sake of simplicity, $r_n=r_n(2\nu-2H,\tauB_{p_k})$ defined by~(\ref{defrnAlphalssgpnu2}) and~(\ref{defrnAlphalssgpnu1}). Note that from Remark~\ref{remhp} $\tauB_{p_k}=2$ for all $k=1,\ldots,K$.

With some little computation, we can obtain, almost surely
\begin{equation}\label{eq1ii}
\! \qtEst{}{p_k}{|Y^{a^m}|^\alpha}  -  \qt{|Y|^\alpha}{p_k} = \alpha \qt{|Y|}{p_k}^{\alpha-1} \left(
\qtEst{}{p_k}{|Y^{a^m}|} - \qt{|Y|}{p_k} 
\right) + \Oas{r_n},
\end{equation}
and
\begin{equation}  \label{eq2ii} 
\qtEst{}{p_k}{\log|Y^{a^m}|} \!\! - \!\!\qt{\log|Y|}{p_k} = \qt{|Y|}{p_k}^{-1} \left(
\qtEst{}{p_k}{|Y^{a^m}|} - \qt{|Y|}{p_k} 
\right) + \Oas{r_n}.
\end{equation}

From~(\ref{eq1-Halpha}), (\ref{eq1-Hlog}), (\ref{eq1ii}), (\ref{eq2ii}) and properties of Gaussian variables, the following results hold almost surely
\begin{equation} \label{eq3ii}
\widehat{H}^{\alpha}-H = \sum_{m=1}^M \sum_{k=1}^K 
\frac{B_m \; c_k}{2 q_k \phi(q_k)} \pi_{k,\alpha} 
\left(
\cdfEst{}{q_k}{|Y|} - p_k \right) + \Oas{r_n} + \OO{b_n},
\end{equation}
and
\begin{equation} \label{eq4ii}
\widehat{H}^{\log}-H = \sum_{m=1}^M \sum_{k=1}^K 
\frac{B_m \; c_k}{2 q_k \phi(q_k)}
\left(
\cdfEst{}{q_k}{|Y|} - p_k \right) + \Oas{r_n} + \OO{b_n},
\end{equation}
where $q_k$ and $\pi_k^\alpha$ are defined by~(\ref{defPikAlpha}). Denote by $\theta_{m,k}^\alpha$ the following constant
$$
\theta_{m,k}^\alpha = \frac{B_m c_k}{2q_k \phi(q_k)} \pi_k^\alpha.
$$
Since $\pi_k^0=1$, (\ref{eq3ii}) and (\ref{eq4ii}) can be rewritten as
\begin{eqnarray} 
\widehat{H}^{\alpha}-H &=& Z_n^\alpha+ \Oas{r_n} + \OO{b_n}
\label{HalphaZn} \\ 
\widehat{H}^{\log}-H &=& Z_n^0 + \Oas{r_n} + \OO{b_n}, \label{HlogZn}
\end{eqnarray}
where for $\alpha \geq 0$, 
\begin{equation} \label{defZnAlpha}
Z_n^\alpha = \sum_{m=1}^M \sum_{k=1}^K \theta_{m,k}^\alpha \left( \cdfEst{}{q_k}{|Y|} - p_k \right) .
\end{equation}
Thus, under Assumption $\bm{A_6(\eta)}$, we have, as $n \to +\infty$,
\begin{eqnarray} 
MSE(\widehat{H}^{\alpha}-H) &=& \OO{\Esp\left((Z_n^\alpha)^2\right)}+ \OO{r_n(2\nu-2H,2)^2} + \OO{n^{-2\eta}},
\label{eq1mse} \\ 
MSE(\widehat{H}^{\log}-H )&=& \OO{\Esp\left((Z_n^0)^2\right)}+ \OO{r_n(2\nu-2H,2)^2} + \OO{n^{-2\eta}}. \label{eq2mse}
\end{eqnarray}
Now,
$$
\Esp\left( (Z_n^\alpha)^2 \right) = \frac{1}{n^2} \sum_{m_1,m_2=1}^M \sum_{k_1,k_2=1}^K \sum_{i_1,i_2=1}^n \theta_{m_1,k_1}^\alpha \theta_{m_2,k_2}^\alpha \Esp\left(
h_{q_{k_1}} ( Y^{\Vect{a}^{m_1}}(i_1) )  h_{q_{k_2}} ( Y^{\Vect{a}^{m_2}}(i_2) )
\right).
$$
For $k_1,k_2=1,\ldots, K$, $m_1,m_2=1,\ldots,M$ and $i_1,i_2=1,\ldots,n$, we have from Lemma~\ref{hp},
\begin{eqnarray}
\Esp\left( 
h_{q_{k_1}} ( Y^{\Vect{a}^{m_1}}(i_1) )  h_{q_{k_2}} ( Y^{\Vect{a}^{m_2}}(i_2) )
\right)& =& \sum_{j_1\geq \tau_{p_{k_1}}/2} \; \sum_{j_2\geq \tau_{p_{k_2}}/2} \frac{c_{2j_1}^{h_{p_{k_1}}} c_{2j_2}^{h_{p_{k_2}}} }{(2j_1)! (2j_2)!} \nonumber \\
&& \qquad \times \Esp \left( 
H_{2j_1}( Y^{\Vect{a}^{m_1}}(i_1) )  H_{2j_2} ( Y^{\Vect{a}^{m_2}}(i_2) ) 
\right) \nonumber \\
&=& \sum_{j \geq 1} \frac{c_{2j}^{h_{p_{k_1}}} c_{2j}^{h_{p_{k_2}}}
 }{(2j)!} \rho^{\Vect{a}^{m_1},\Vect{a}^{m_2}} (i_2-i_1)^{2j}. \label{eq-calculhqCorr}
\end{eqnarray}
Under Assumption $\bm{A_1(2\nu)}$, we have from Lemma~\ref{lemKent}, $\rho^{\Vect{a}^{m_1},\Vect{a}^{m_2}} (i)=\OO{|i|^{2H-2\nu}}$. Now, we leave the reader to check that, as $n \to +\infty$
$$
\frac{1}{n^2} \sum_{i_1,i_2=1}^n \rho^{\Vect{a}^{m_1},\Vect{a}^{m_2}} (i_2-i_1)^{2}
= \OO{\frac1n \sum_{|i|\leq n} |i|^{2(2H-2\nu)}} = \OO{v_n(2\nu-2H))},
$$
where the sequence $v_n(\cdot)$ is given by~(\ref{defvn}). Thus, we have, as $n \to +\infty$, $\Esp\left( (Z_n^\alpha)^2 \right)= \OO{v_n(2\nu-2H)}$, which leads to the result from~(\ref{eq1mse}) and~(\ref{eq2mse}).

$(iii)$ Assume $\nu>H+1/4$ and $\eta>1/2$, then from~(\ref{HalphaZn}) and~(\ref{HlogZn}), the following equivalences in distribution hold
\begin{equation} \label{equivZnHest}
\sqrt{n} \left( \widehat{H}_n^{\alpha}-H \right) \sim \sqrt{n} \; Z_n^\alpha \quad \mbox{ and } \quad
\sqrt{n} \left( \widehat{H}_n^{\log}-H \right) \sim \sqrt{n} \; Z_n^0.
\end{equation}
Now, decompose $Z_n^\alpha=T_n^1 + T_n^2$, where  
$$
T_n^{1} = \frac1{\sqrt{n}} \sum_{m=1}^M \sum_{k=1}^K \theta_{m,k}^\alpha  \sum_{i=\ell+1}^{M\ell+1} h_{q_k}(Y^{\Vect{a^m}}(i)) $$
and 
$$T_n^{2} =  \sqrt{n} \sum_{m=1}^M \sum_{k=1}^K \theta_{m,k}^\alpha \left\{\frac1n \sum_{i=M\ell+1}^n h_{q_k}(Y^{\Vect{a^m}}(i)) \right\},
$$
Clearly, $T_n^{1}$ converges to 0 in probability, as $n \to +\infty$. Therefore, we have, as $n\to +\infty$
\begin{equation} \label{Tn2Alpha}
Z_n^\alpha \sim \sqrt{n} \left\{  \frac1n \sum_{i=M\ell+1}^n G^\alpha \left( Y^{\Vect{a^1}}(i),\ldots, Y^{\Vect{a^M}}(i)\right) \right\}
\end{equation}
where $G^\alpha$ is the function from $\RR^M$ to $\RR$ defined for $\alpha\geq 0$ and $t_1,\ldots,t_M \in \RR$ by: 
\begin{equation} \label{defGAlpha}
G^\alpha(t_1,\ldots,t_M)=\sum_{m=1}^M\sum_{k=1}^K \theta_{m,k} \; h_{q_k}(t_m).
\end{equation}
Denote by $\Vect{\widetilde{Y}^a}(i)$, the vector defined for $i=M\ell+1,\ldots,n$ by 
$$\Vect{\widetilde{Y}^a}(i)= (  Y^{\Vect{a^1}}(i),\ldots, Y^{\Vect{a^M}}(i) ).$$
We obviously have $\Esp\big(G^\alpha (\Vect{\widetilde{Y}^a}(i))^2\big)< + \infty$. Since, for all $k=1,\ldots,K$, the functions $h_{q_k}$ have Hermite rank~$\tau_{p_k}$, the function $G^\alpha$ has Hermite rank $2$ (see {\it e.g.} \Citet{Arcones94} for the definition of the Hermite rank of multivariate functions). Moreover under Assumption $\bm{A_1(2\nu)}$, we have from Lemma~\ref{lemKent}, as $j \to +\infty$
\begin{eqnarray}
\Esp \left( Y^{\Vect{a^{m_1}}}(i) Y^{\Vect{a^{m_2}}}(i+j) \right) &=&
\rho^{\Vect{ a^{m_1},a^{m_2}}} (j)
= \OO{|j|^{2H-2\nu}} \in \ell^2(\ZZ),  \nonumber
\end{eqnarray}
as soon as  $\nu > H+1/4$. Thus, from Theorem~4 of \Citet{Arcones94}, there exists $\sigma^2_\alpha$  (defined for $\alpha \geq 0$) such that, as $n\to +\infty$, the following convergence in distribution 
holds
$$
Z_n^\alpha  \longrightarrow \mathcal{N}(0,\sigma^2_{\alpha}) $$
with 
$$\sigma^2_{\alpha} = \sum_{i \in \ZZ} \Esp \left( G^\alpha \left(  \Vect{\widetilde{Y}^a}(i^\prime) \right)  G^\alpha \left(  \Vect{\widetilde{Y}^a}(i^\prime+i) \right) \right).$$
With previous notations, we have
\begin{eqnarray}
\sigma^2_{\alpha} &=& \sum_{i\in\ZZ} \sum_{m_1,m_2=1}^{M} \sum_{k_1,k_2=1}^K \theta_{m_1,k_1}^\alpha \theta_{m_2,k_2}^\alpha  
\Esp \left( h_{p_{k_1}}( Y^{\Vect{a^{m_1}}}(i^\prime) )  h_{p_{k_2}}( Y^{\Vect{a^{m_2}}}(i^\prime+i) )     \right) \nonumber \\
%&=& \sum_{i\in\ZZ} \sum_{m_1,m_2}^{M} \sum_{k_1,k_2}^K  \sum_{j_1,j_2 \geq r} c_{2j_1}^{h_{p_{k_1}}} c_{2j_2}^{h_{p_{k_2}}}  \lambda_{m_1,k_1}^\alpha \lambda_{m_2,k_2}^\alpha  
%\Esp \left(  H_{2j_1} (  Y^{\Vect{a^{m_1}}}(i^\prime) )    H_{2j_2} (  Y^{\Vect{a^{m_2}}}(i^\prime+i) )        \right) \nonumber \\
&=& \sum_{i\in\ZZ} \sum_{m_1,m_2=1}^{M} \sum_{k_1,k_2=1}^K  \sum_{j\geq r} \frac{c_{2j}^{h_{p_{k_1}}} c_{2j}^{h_{p_{k_2}}} }{(2j)!}  \theta_{m_1,k_1}^\alpha \theta_{m_2,k_2}^\alpha \; \rho^{\Vect{ a^{m_1},a^{m_2}}} (i)^{2j}. \label{lastExpr}
\end{eqnarray}
From (\ref{coeffhp}), we can see that formula (\ref{lastExpr}) is equivalent to (\ref{sig2Alpha}), which ends the proof from~(\ref{equivZnHest}).
\end{proof}

\subsection{Proof of Corollary~\ref{corolHest}}

\begin{proof}
Equation (\ref{Tn2Alpha}) is still available for a sequence $\alpha_n$ such that $\alpha_n \to 0$ as $n \to +\infty$, that is 
$$
\sqrt{n} \left( \widehat{H}_n^{\alpha_n} - H \right) \sim \sqrt{n} \left\{  \frac1n \sum_{i=M\ell+1}^n G^{\alpha_n} \left( Y^{\Vect{a^1}}(i),\ldots, Y^{\Vect{a^M}}(i)\right) \right\}
$$
From (\ref{defGAlpha}) and since $\pi_k^{\alpha_n} \to 1$, as $n \to +\infty$, we have $G^{\alpha_n}(\cdot) \to G^{0}(\cdot)$. Therefore, the following equivalence in distribution holds, as $n \to +\infty$
$$
\sqrt{n} \left( \widehat{H}_n^{\alpha_n} - H \right) \sim \sqrt{n} \left( \widehat{H}_n^{\log} - H \right),
$$
which ends the proof.
\end{proof}

\subsection{Auxiliary Lemma for the proof of Theorem~\ref{convHestTM}}

\begin{lemma} \label{lemme-auxTM}
Let $0<\beta_1\leq\beta_2<1$ and let $\Vect{Z}=(Z_1,\ldots,Z_n)$ $n$ random variables identically distributed, such that
$\sup_{\beta_1\leq p \leq \beta_2}\qtEst{Z}{p}{Z}=\Oas{1}$, then
$$
\tmEst{Z} - \frac1{1-\beta_2-\beta_1}\int_{\beta_1}^{1-\beta_2}\qtEst{Z}{p}{Z}dp = \Oas{n^{-1}}.
$$
\end{lemma}

\begin{proof}
It is sufficient to notice that for $i=1,\ldots,n$
$$
n \int_{\frac{i-1}n}^{\frac{i}n} \qtEst{}{p}{Z}dp \leq Z_{(i),n} \leq n \int_{\frac{i}n}^{\frac{i+1}n} \qtEst{}{p}{Z}dp.
$$
which leads to
$$
\frac{n}{n-[n\beta_2]-[n\beta_1]} \int_{\frac{[n\beta_1]}n}^{\frac{n-[n\beta_2]}{n}} \qtEst{}{p}{Z} dp\leq \tmEst{Z}\leq 
\frac{n}{n-[n\beta_2]-[n\beta_1]} \int_{\frac{[n\beta_1]}n+\frac1n}^{\frac{n-[n\beta_2]}{n}+\frac1n} \qtEst{}{p}{Z} dp.
$$
The end is omitted.
\end{proof}

\subsection{Proof of Theorem~\ref{convHestTM}}

\begin{proof}
$(i)$ From~(\ref{defepsmalphaTM}), (\ref{defepsmlogTM}), and~(\ref{HestMoinsHTM}), we have 
\begin{eqnarray}
\hspace*{-0.2cm}\widehat{H}^{\alpha,tm}\!\!-\!\!H &=& \sum_{m=1}^M \frac{B_m}{\alpha} \varepsilon_m^{\alpha,tm}  \nonumber \\
&=&\sum_{m=1}^M \frac{B_m}{\alpha} \log \left( \tmEst{|Y^{a^m}|^\alpha}/ \tm{|Y|^\alpha} \right)  +\alpha \times b_n \nonumber \\
&=& \sum_{m=1}^M \frac{B_m}{\alpha \tm{|Y|^\alpha}} \left( \tmEst{|Y^{a^m}|^\alpha}- \tm{|Y|^\alpha}\right) \!(1+\oas{1}) + \alpha \; b_n. \label{eq1-HalphaTM} 
\end{eqnarray}
and
\begin{eqnarray}
\widehat{H}^{\log,tm}-H &=& \sum_{m=1}^M  B_m \varepsilon_m^{\log,tm} \nonumber \\
&=& \sum_{m=1}^M B_m \left( 
\tmEst{\log|Y^{a^m}|}-\tm{\log|Y|} \right) 
+ b_n \label{eq1-HlogTM}
\end{eqnarray}
Let us notice that from~Lemma~\ref{lemme-supqQuant}, one can apply Lemma~\ref{lemme-auxTM} for the vectors $\Vect{|Y^{a^m}|^\alpha}$ and $\Vect{\log |Y^{a^m}|}$. Then it comes
\begin{eqnarray*}
\tmEst{|Y^{a^m}|^\alpha}\!\! -\!\!\tm{|Y|^\alpha} &=&\!\! \frac1{1-\beta_2-\beta_1} \int_{\beta_1}^{1-\beta_2} \left(\qtEst{}{p}{|Y^{a^m}|^\alpha} - \qt{|Y|^\alpha}{p}  \right)dp + \Oas{n^{-1}},  \\
\!\!\!\tmEst{\log |Y^{a^m}|}\!\! -\!\! \tm{\log |Y|} &=&\!\! \frac1{1-\beta_2-\beta_1} \int_{\beta_1}^{1-\beta_2} \!\!\!\left(\qtEst{}{p}{\log |Y^{a^m}|} - \qt{\log |Y|}{p}  \right)dp + \Oas{n^{-1}}.
\end{eqnarray*}
Hence, from (\ref{rappelCorrelation}) and Lemma~\ref{lemme-supqQuant} and Remark~\ref{remhp} and under Assumption~$\bm{(A_6(\eta))}$, we obtain
\begin{eqnarray*}
\widehat{H}^{\alpha,tm}-H &=& \Oas{y_n(2H-2,2)} + \OO{n^{-\eta}} + \Oas{n^{-1}} \\ \widehat{H}^{\log,tm}-H &=&  \Oas{y_n(2H-2,2)} + \OO{n^{-\eta}} + \Oas{n^{-1}},
\end{eqnarray*}
where the sequence $y_n(\cdot,\cdot)$ is defined by~(\ref{defeps}) with $L(\cdot)=1$. This leads to the result by noticing that $n^{-1}=\OO{y_n({2H-2},2)}$.

$(ii)$ By following the proof of Theorem~\ref{convHest}~$(ii)$ and from Theorem~\ref{unifBahadur}, we may obtain the following representation
\begin{eqnarray*}
\widehat{H}^{\alpha,tm}-H &=& \sum_{m=1}^M \frac{B_m}{\alpha \tm{|Y|^\alpha}} \times \frac1{1-\beta_2-\beta_1} \int_{\beta_1}^{1-\beta_2} \frac{
\cdfEst{}{q}{|Y^{a^m}|} - p}{2 \frac{1}{\alpha} q^{1-\alpha} \phi(q)} dp \\
&& + \Oas{r_n}+ \Oas{n^{-1}} + \OO{n^{-\eta}}, \\
\widehat{H}^{\log,tm}-H &=& \sum_{m=1}^M B_m \times \frac1{1-\beta_2-\beta_1} \int_{\beta_1}^{1-\beta_2} \frac{
\cdfEst{}{q}{|Y^{a^m}|} - p}{2 q \phi(q)} dp \\
&&+ \Oas{r_n}+ \Oas{n^{-1}} + \OO{n^{-\eta}}. \\
\end{eqnarray*}
With such a representation, we observe that the result $(ii)$ can be proved similarly to the one of Theorem~\ref{convHest}. 

$(iii)$ By assuming that $\eta>1/2$ and $\nu>H+1/4$, one may obtain the asymptotic normality of $\widehat{H}^{\alpha,tm}$ and $\widehat{H}^{\log,tm}$ by using the same tools as the one presented in the proof of Theorem~\ref{convHest}~$(iii)$. Therefore, let us just explicit the asymptotic variance of estimators $\widehat{H}^{\alpha,tm}$ and $\widehat{H}^{\log,tm}$.
If $\nu>H+1/4$ and $\eta>1/2$, then from previous representations and from~\ref{eq-calculhqCorr} we obtain as $n \to +\infty$
\begin{eqnarray*}
Var\left( \sqrt{n} \left( \widehat{H}^{\alpha,tm}-H\right)\right) &\sim& \frac n{n^2} \sum_{m_1,m_2=1}^M \sum_{i_1,i_2=1}^n  B_{m_1}B_{m_2} \frac1{(\tm{|Y|^\alpha})^2}\times \frac1{(1-\beta_2-\beta_1)^2} \times \\ && \int_{\beta_1}^{1-\beta_2}  \int_{\beta_1}^{1-\beta_2} \sum_{j\geq 1} \frac{c_{2j}^{h_{p_{1}}} c_{2j}^{h_{p_{2}}} }{(2j)!} \frac{q_1^{\alpha-1}q_2^{\alpha-1}}{4 \phi(q_1) \phi(q_2)}dp_1dp_2 \rho^{a^{m_1},a^{m_2}}(i)^{2j} ,
\end{eqnarray*}
with $q_k=\Phi^{-1}\left(\frac{1+p_k}{2}\right)$ for $k=1,2$. Due to~(\ref{coeffhp}) and since
$\tm{|Y|^\alpha}=\frac1{1-\beta_2-\beta_1}\int_{\beta_1}^{1-\beta_2}q^\alpha dp$, this variance converges towards $\sigma_{\alpha,tm}^2$ given by~(\ref{sigma2TM}), as $n \to +\infty$. 

We leave the reader to check that the asymptotic variance of $\sqrt{n} \left( \widehat{H}^{\log,tm}-H\right)$ is given by $\sigma^2_{0,tm}$.
\end{proof}

%\begin{center}
%\begin{tabular}{ll}
%\hline
%Type Robustesse & Effet \\ \hline
%0 & processus non contaminé \\
%7 & deux observations sont perturbées ($\sigma=3*\sigma_{proc}$) \\
%9& deux observations sont perturbées ($\sigma=6*\sigma_{proc}$) \\
%10 &  quatre observations sont perturbées ($\sigma=6*\sigma_{proc}$)\\
%11 & processus bruité ($\sigma=0.02$)\\
%\hline
%\end{tabular}
%\end{center}

%\input textFGNTypes
%\input textEXPTypes
%\input textLOGTypes
%\input textARFTypes

\bigskip

{\bf \large Acknowledgement.} The author is very grateful to Anestis Antoniadis and R\'emy Drouilhet for helpful comments and to Kinga Sipos for a careful reading of the present paper.

\bigskip

\noindent {\sc J.-F. Coeurjolly \\
LJK, SAGAG Team, Universit\'e Grenoble 2 \\
1251 Av. Centrale  BP 47\\
38040 GRENOBLE Cedex 09\\
France \\
E-mail:} Jean-Francois.Coeurjolly@upmf-grenoble.fr

%\newpage 

\begin{table}[ht]
\begin{center}
\begin{tabular}{lcc}
\hline\hline
\multicolumn{3}{c}{Non-contaminated sample paths}\\
\hline
Estimators & {\scriptsize$v(\cdot)=|\cdot|^{2H}$ }& {\scriptsize$v(\cdot)=1-\exp(-|\cdot|^{2H})$} \\
\hline
$\Vect{p}=1/2$, $\Vect{c}=1$ (median)& 0.796 (0.042) & 0.801 (0.042) \\
$\Vect{p}=0.9$, $\Vect{c}=1$ & 0.797 (0.035) & 0.798 (0.036) \\
$\Vect{p}=(1/4,3/4)$,$\Vect{c}=(1/2,1/2)$, $g(\cdot)=|\cdot|^2$ & 0.795 (0.036) & 0.800 (0.037) \\
%$\Vect{p}=(1/4,3/4)$,$\Vect{c}=(1/2,1/2)$ , $g(\cdot)=\log|\cdot|$   & 0.797 (0.042) & 0.796 (0.032) \\
10\%-trimmed mean, $g(\cdot)=|\cdot|^2$ & 0.797 (0.03) & 0.799 (0.034) \\
Quadratic variations method& 0.802 (0.032) & 0.798 (0.032) \\
Whittle estimator & 0.805 (0.024) & 0.806 (0.024) \\
\hline\hline
\end{tabular}

\bigskip

\begin{tabular}{lcc}
\hline\hline
\multicolumn{3}{c}{Contaminated sample paths}\\
\hline
Estimators & {\scriptsize$v(\cdot)=|\cdot|^{2H}$ }& {\scriptsize$v(\cdot)=1-\exp(-|\cdot|^{2H})$} \\
\hline
$\Vect{p}=1/2$, $\Vect{c}=1$ (median)& 0.798 (0.047) & 0.803 (0.045) \\
$\Vect{p}=0.9$, $\Vect{c}=1$ & 0.793 (0.033) & 0.789 (0.032) \\
$\Vect{p}=(1/4,3/4)$,$\Vect{c}=(1/2,1/2)$, $g(\cdot)=|\cdot|^2$ & 0.797 (0.040) & 0.796 (0.037) \\
%$\Vect{p}=(1/4,3/4)$,$\Vect{c}=(1/2,1/2)$ , $g(\cdot)=\log|\cdot|$   & 0.798 (0.044) & 0.804 (0.044) \\
10\%-trimmed mean, $g(\cdot)=|\cdot|^2$ & 0.792 (0.037) & 0.797 (0.033) \\
Quadratic variations method& 0.329 (0.162) & 0.353 (0.149) \\
Whittle estimator & 0.519 (0.106) & 0.510 (0.100) \\
\hline\hline
\end{tabular}
\end{center}
\caption{\small Mean and standard deviations for $n=1000$ and $H=0.8$ using $500$ Monte Carlo simulations of sample paths of processes with variance function $v(\cdot)=|\cdot|^{2H}$, respectively $v(\cdot)=1-\exp(-|\cdot|^{2H})$ (first table ) and contaminated versions (second table), see (\ref{contamination}).}
\label{tab-Robuste}
\end{table}

\end{document}